\documentclass[11pt]{article}

\usepackage{amsmath,amsthm,amsfonts,latexsym}
\usepackage{amssymb}
\usepackage[latin1]{inputenc}
\usepackage{tikz-cd}

\begin{document}

\newtheorem*{theo}{Theorem}
\newtheorem*{pro} {Proposition}
\newtheorem*{cor} {Corollary}
\newtheorem*{lem} {Lemma}
\newtheorem{theorem}{Theorem}[section]
\newtheorem{corollary}[theorem]{Corollary}
\newtheorem{lemma}[theorem]{Lemma}
\newtheorem{proposition}[theorem]{Proposition}
\newtheorem{conjecture}[theorem]{Conjecture}

\theoremstyle{definition}
 \newtheorem{definition}[theorem]{Definition} 
  \newtheorem{example}[theorem]{Example}
   \newtheorem{remark}[theorem]{Remark}
   
\newcommand{\Naturali}{{\mathbb{N}}}
\newcommand{\Reali}{{\mathbb{R}}}
\newcommand{\Complessi}{{\mathbb{C}}}
\newcommand{\Toro}{{\mathbb{T}}}
\newcommand{\Relativi}{{\mathbb{Z}}}
\newcommand{\HH}{\mathfrak H}
\newcommand{\KK}{\mathfrak K}
\newcommand{\LL}{\mathfrak L}
\newcommand{\as}{\ast_{\sigma}}
\newcommand{\tn}{\vert\hspace{-.3mm}\vert\hspace{-.3mm}\vert}
\def\A{{\mathcal A}}
\def\B{{\mathcal B}}
\def\E{{\mathcal E}}
\def\F{{\mathcal F}}
\def\H{{\mathcal H}}
\def\K{{\mathcal K}}
\def\L{{\mathcal L}}
\def\N{{\mathcal N}}
\def\M{{\cal M}}
\def\gM{{\frak M}}
\def\O{{\cal O}}
\def\P{{\cal P}}
\def\S{{\cal S}}
\def\T{{\cal T}}
\def\U{{\cal U}}
\def\C{{\mathcal C}}
\def\V{{\mathcal V}}
\def\X{{\mathcal X}}
\def\Y{{\mathcal Y}}
\def\qed{\hfill$\square$}
\def\veps{{\varepsilon}}

\title{Positive definiteness and Fell bundles over discrete groups}

\author{Erik B\'edos, Roberto Conti\\}
\date{\today}
\maketitle
\markboth{R. Conti, Erik B\'edos}{
}
\renewcommand{\sectionmark}[1]{}
\begin{abstract} We introduce a natural concept of positive definiteness for bundle maps between Fell bundles over (possibly different) discrete groups and describe several examples. Such maps induce completely positive maps between the associated full cross-sectional C*-algebras in a functorial way. Under the assumption that the kernel of the homomorphism connecting the groups under consideration is amenable, they also induce completely positive maps between the associated reduced cross-sectional C*-algebras. As an application, we define an approximation property for a Fell bundle over a discrete group which generalizes Exel's approximation property and still implies the weak containment property. Both approximation properties coincide when the unit fibre is nuclear.
\end{abstract}



\section{Introduction}
  The concept of positive definiteness for complex-valued functions on a group $G$ has a long history 
  (it goes back to Toeplitz in 1911 when $G=\Relativi$ and to Mathias in 1923 when $G=\Reali$; for more details, see \cite{Ste} and references therein). Its importance is ubiquitous in various areas of mathematics, in particular in connection with the development in the 1940's and 50's of the representation theory of locally compact groups and of the associated group $C^*$-algebras  (see e.g.~\cite{Dix, FD1, FD2}). On the other hand, completely positive linear maps play a prominent r\^ole in the modern theory of operator algebras and in quantum information theory (cf.~e.g.~\cite{Pau, BrOz, Hay}). Haagerup discovered in 1978 that there is a very useful link between the two concepts. Namely, he proved in \cite[Theorem 3.1]{Ha1} that given an action $\alpha$ of a locally compact group $G$ on a von Neumann algebra $\M$, every continuous positive definite function $\psi$ on $G$ gives rise to a normal completely positive linear map $M_\psi$ on the crossed product  $\M\rtimes_\alpha G$ such that $M_\psi$ acts on each $g$-component by multiplication with $\psi(g)$. This link also holds in the context of group $C^*$-algebras and group von Neumann algebras, as noted by Haagerup in his seminal paper \cite{Ha2} and exploited in later works, such as \cite{DCHa} where continuous positive definite functions on $G$ are identified as multipliers of the Fourier algebra $A(G)$ of
  $G$.
  In the setting of discrete unital $C^*$-dynamical systems, a similar link was established in a series of paper \cite{BeCo3, BeCo4, BeCo6}, where broader concepts of multipliers were introduced. As a sample of related articles we mention \cite{McKTT, McKSTT, McKT, HTT, BC, KLS, McKPTT}.  

Given two $C^*$-algebras $A$ and $B$, it is often the case that
one needs to construct completely positive linear maps from $A$ to $B$. 
When $A$ and $B$ belong to some class of $C^*$-algebras, one may hope to find some specific procedure achieving this. In this paper, we consider this problem for the class $\mathcal{F}$ (resp.~$\mathcal{F}_{\rm r}$) of $C^*$-algebras which may be described as the full (resp.~reduced) cross-sectional $C^*$-algebra $C^*(\A)$ (resp.~$C_r^*(\A)$) of a Fell bundle $\A$ over a discrete group, as studied in \cite{Exel1, Exel2} (and e.g.~in \cite{KS, KM}; see also \cite{FD1, FD2} for a more general treatment dealing with $C^*$-algebraic bundles, i.e., Fell bundles over locally compact groups). The class 
$\mathcal{F}$ (resp.~$\mathcal{F}_{\rm r}$) is large as it contains all  full (resp.~reduced) twisted group $C^*$-algebras, and more generally all full (resp.~reduced) $C^*$-crossed products associated to discrete twisted $C^*$-dynamical systems, even those for which the action of the group is only assumed to be partial \cite{Exel0}. Also, every  $C^*$-algebra on which there exists a maximal (resp.~normal) nondegenerate coaction by a discrete group belongs to $\mathcal{F}$ (resp.~$\mathcal{F}_{\rm r}$), see \cite[Proposition 4.2]{EKQ} (resp.~\cite[Corollary 3.9]{Q}). In particular,  both classes include any $C^*$-algebra on which there exists a nondegenerate coaction by a discrete amenable group,  e.g., any $C^*$-algebra on which there exists an action of a compact abelian group, such as
Cuntz algebras 
and, more generally, 
(higher-rank) graph $C^*$-algebras. Note here that one may always choose to regard 
a $C^*$-algebra $B$ as the cross-sectional $C^*$-algebra of a Fell bundle over the trivial group, hence that our set up includes the case where no assumption is put on $B$. 
 
     After a preliminary section (Section 2), our main findings are gathered in Section 3. To describe these, let us first assume that $A= C^*(\A)$ and $B=C^*(\B)$, where $\A$ is a Fell bundle over a discrete group $G$ and $\B$ is a Fell bundle over a discrete group $H$, and let $\varphi:G\to H$ be a homomorphism. Consider an $\A$-$\varphi$-$\B$ bundle map $T:\A\to \B$, meaning that $T$ maps the fibre of $\A$ over any $g\in G$ into the fibre of  $\B$ over $\varphi(g)$ in a linear and bounded way. Inspired by the definition of positive definiteness for a complex function on a group and the definition of complete positivity for a linear map between $C^*$-algebras, we introduce a natural notion of positive definiteness for such a bundle map, and show that $T$ is positive definite if and only if there exists a completely positive linear map $\Phi_T:A\to B$ which ``extends'' $T$ (see Theorem \ref{cp-full-bundle} for a more precise statement). Our proof relies on a reformulation of positive definiteness in terms of the matrix $C^*$ algebras associated to a Fell bundle by Abadie and Ferraro in \cite{AF}, and builds on the Stinespring dilation theorem for Fell bundles due to Buss, Ferraro and Sehnem in \cite[Appendix]{BFS}. 
     As we point out in Remark \ref{category}, the association $\A\mapsto C^*(\A), \,(\varphi, T) \mapsto \Phi_T$ is functorial. 
  
 Next, we consider the case where $A= C_r^*(\A)$ and $B=C_r^*(\B)$. It is not difficult to see that an $\A$-$\varphi$-$\B$ bundle map $T$ has to be positive definite whenever there exists a completely positive linear map $M_T:A\to B$ which ``extends'' $T$. However, the converse statement can not hold in general. This may be deduced from a result in the recent book of Bekka and de la Harpe, cf.~\cite[Corollary 8.C.15]{BedlH}, which says that the homomorphism $\varphi:G\to H$ extends to a $*$-homomorphism from $C_r^*(G)$ into $C_r^*(H)$ if and only if the kernel of $\varphi$ is an amenable subgroup of $G$. Using this result in combination with Fell's absorption principle for Fell bundles and a generalization of another result of Buss, Ferraro and Sehnem (\cite[Proposition 4.8]{BFS}), we are able to show that the converse does hold if we also assume that $\ker(\varphi)$ is amenable, cf.~Theorem \ref{cp-reduced}. Again, the association $\A\mapsto C_r^*(\A),  (\varphi, T) \mapsto M_T$ is functorial in a suitable sense.

 In Section 4 we illustrate the usefulness of Theorem \ref{cp-full-bundle} and Theorem \ref{cp-reduced} in a series of examples, showing that these two results are not only generalizations of previous results of the same flavor, but also lead to new applications. 
 In particular, we devote a section (Section 5) to discussing a new approximation property for a Fell bundle $\A$ over a discrete group. We call it the PD-approximation property since it requires the existence of a net of uniformly bounded, finitely supported positive definite bundle maps converging to the identity map on each fibre, in analogy with the well-known characterization of the amenability of a discrete group in terms of positive definite functions. The  PD-approximation property for $\A$ is (formally) weaker than Exel's approximation property \cite{Exel1, Exel2}, but still implies the weak containment property (i.e., amenability of $\A$ in the sense of Exel), cf.~Theorem \ref{amen}. Moreover, when it is satisfied, then  $C^*(\A)$ (resp.~$C^*_r(\A)$) is nuclear if and only if the unit fibre of $\A$ is nuclear, cf.~Theorem \ref{amen2}. Using current knowledge from \cite{ABF2}, it follows that the PD-approximation property for $\A$ is equivalent to Exel's approximation property when the unit fibre of $\A$ is assumed to be nuclear. However, the arguments used to deduce 
this fact
 give no clue whether it also holds
 when the unit fibre of $\A$ is not nuclear.

 In the final section (Section 6) we discuss some possible further developments. Let us mention here another one. 
  Abadie and Ferraro have introduced in \cite{AF} (see also \cite{ABF1}) the notion of a right Hilbert $\B$-bundle for a Fell bundle $\B$ over a group, which reduces to the notion a right Hilbert $C^*$-module when the group is trivial.  In a forthcoming article we will consider left actions of Fell bundles on such right Hilbert bundles and discuss the connection with positive definite bundle maps. This will provide a useful tool to produce positive definite bundle maps.
  Moreover, we will also study how $C^*$-correspondences over the cross-sectional $C^*$-algebras arise from such actions.

When we were about to finish the first draft
of the present article, 
an interesting preprint by Buss, Kwasniewski, McKee and Skalski appeared  \cite{BKMS}, where the existing theories of Fourier-Stieltjes algebras associated to groups \cite{KL}, to groupoids \cite{Re, Oty, Pat} and to discrete twisted unital $C^*$-dynamical systems \cite{BeCo6}  is extended to the category of twisted actions by \'etale groupoids on $C^*$-bundles. Only a few of their results deal with the general case of Fell bundles over \'etale groupoids, hence with general Fell bundles over discrete groups, and 
there is little overlap 
between our papers.

\section{Preliminaries} \label{prem}
\medskip Throughout this article, $G$ and $H$ will always denote discrete groups. The unit of any group will be denoted by $e$. We will be considering Fell bundles over discrete groups. Our main reference on this topic will be Exel's book \cite{Exel2}, and we will follow his notation and terminology. For the ease of the reader we give below a short introduction.  We assume familiarity with Hilbert $C^*$-modules, and follow the notation in Lance's book \cite{La1}. In particular, if $X$ is a right Hilbert module over a $C^*$-algebra $A$, $\L_A(X)$ will denote the $C^*$-algebra of adjointable operators on $X$. The left regular representation of $G$ on $\ell^2(G)$ will be denoted by $\lambda^G$.

\begin{definition}(\cite[Definition 16.1]{Exel2})
A \emph{Fell bundle $\A=(A_g)_{g\in G}$ over $G$} is a collection of Banach spaces (called \emph{fibers}) satisfying the following properties (where $\A$ also denotes the disjoint union of the $A_g$'s).
There is an associative multiplication map from $\A\times \A$ into $\A$ 
such that $A_g A_h \subseteq A_{gh}$, $(a, b)\mapsto ab$ is bilinear on $A_g\times A_h$ for all $g, h \in G$, and $\|ab\| \leq \|a\|\|b\|$ for all $a, b\in A$. Moreover, there is an involutive, anti-multiplicative map from $\A$ into itself
such that $A_g^*= A_{g^{-1}}$ and $a\mapsto a^*$ is conjugate-linear, norm-preserving from $A_g$ into $A_{g^{-1}}$ for all $g \in G$.  Finally, we have $\|a^*a\| =\|a\|^2$ and $a^*a \geq 0$ in the unit fibre $A_e$ (which is a $C^*$-algebra) for all $a\in \A$.
\end{definition}
The vector space $C_c(\A)$ consisting of all finitely supported functions $f:G\to \A$ such that $ f(g) \in A_g \text{ for all } g \in G$,
becomes a $*$-algebra
with respect to the operations defined by
\[ (f_1 \ast f_2) (h) = \sum_{g\in G} f_1(g) f_2(g^{-1}h), \quad f^*(h) = f(h^{-1})^*\]
for all $h \in G$. 
 The space $C_c(\A)$ is also a right pre-Hilbert $A_e$-module,  the right action of $A_e$ being defined pointwise and the inner product being given by  \[ \langle f_1, f_2 \rangle_{A_e} = \sum_{h\in G} f_1(h)^*f_2(h).\]
 By completion we obtain a right Hilbert $A_e$-module, denoted by $\ell^2(\A)$.
The (left) regular representation $\lambda^\A = (\lambda^\A_g)_{g\in G}$ of $\A$ in  $\L_{A_e}(\ell^2(\A))$ is determined  for each $g\in G$
by
\[ (\lambda^\A_g(a)f)(h) = a\,f(g^{-1}h) \quad \text{whenever } a \in A_g, f \in C_c(\A) \text{ and } h \in G.\]
It induces a faithful $*$-representation $\iota^\A: C_c(\A)\to \L_{A_e}(\ell^2(\A))$ given  by
\[ \iota^\A(f) = \sum_{g\in G} \lambda^\A_g(f(g)) \quad \text{ for all } f\in C_c(\A),\]
see  \cite[Proposition 17.9, (ii)]{Exel2}.
 
The reduced cross-sectional $C^*$-algebra $C_r^*(\A)$ is the $C^*$-subalgebra of $\L_{A_e}(\ell^2(\A))$ generated by $\iota^\A(C_c(\A))$, i.e., by $\{ \lambda_g^\A(a): g\in G, a \in A_g\}$. 
The full cross-sectional $C^*$-algebra $C^*(\A)$ is defined as the $C^*$-completion of $C_c(\A)$ with respect to the universal $C^*$-norm  $\|\cdot\|_{\rm u}$ given by \[\|f\|_{\rm u} 
=  \sup\big\{p(f): p \text{ is a } C^*\text{-seminorm on } C_c(\A)\big\}.\]
  For $f\in C_c(\A)$, we set $\|f\|_r := \|\iota^\A(f)\|$. We then have $\|f\|_{r} \leq \|f\|_{\rm u}$.

For $g \in G$, we let $j^\A_g: A_g \to C_c(\A)$ denote the canonical injection. If $a\in A_g$ we will sometimes write $a\odot g$ instead of $j^\A_g(a)$. Further, we let $\kappa^\A:C_c(\A)\to C^*(\A)$ 
denote the canonical injection and set $\widehat j^\A_g := \kappa^\A\circ j^\A_g : A_g \to C^*(\A)$. 
   For any $f \in C_c(\A)$, we have $f = \sum_{g\in G} j^{\A}_g(f(g))$,
 so 
\[\kappa^\A(f) =  \sum_{g\in G} \widehat j^\A_g(f(g)).\] 
As shown in \cite[Proposition 17.9, (iv)]{Exel2}, the map $\widehat j^{\A}_g$ is isometric for every $g\in G$.
Thus we get that  
\begin{equation}\label{norm-kappa} \| \kappa^\A(f)\|_{\rm u} \leq 
\sum_{g\in G} \|\widehat{j}^\A_g(f(g))\|_{\rm u}
= \sum_{g\in G} \|f(g)\|.
\end{equation}
The canonical $*$-homomorphism $\Lambda^\A$ from $C^*(\A)$ onto $C_r^*(\A)$ is determined by   \[\Lambda^\A\big(\,\widehat{j}_g^\A(a)\big) = \lambda_g^\A(a)\] for all $g\in G$ and $a \in A_g$, i.e., it satisfies $\iota^\A = \Lambda^\A \circ \kappa^\A$. 
The following diagram allows to visualize 
the various maps:
\[
\begin{tikzcd}
& & C^*(\A) \arrow [dd, two heads ,"\Lambda^\A"]  \\
A_g \arrow [r, "j^\A_g"] \arrow[rru, bend left, "\widehat{j}^\A_g"]  \arrow[rrd, bend right, "\lambda^\A_g"] & C_c(\A) \arrow [ru, "\kappa^\A"] \arrow[rd, "\iota^\A"] \\
& & C^*_r(\A) 
\end{tikzcd}
\]
In addition, for each $g\in G$, there is a contractive linear map $E_g^\A : C_r^*(\A) \to A_g$ satisfying that for each $h\in G$ and $a\in A_h$, we have
\[E_g^\A\big(\lambda^\A_h (a)\big)=\begin{cases} a, \quad \text{ if } g =h, \\  0, \quad \text{ if} g\neq h, \end{cases}\]
cf.~\cite[Lemma 17.8]{Exel2}. For $x \in C_r^*(\A)$, $E_g^\A(x)$ may be thought of as the Fourier coefficient of $x$ at $g$. The map $E^\A := E_e^\A$ is a faithful conditional expectation from $C_r^*(\A)$ onto $A_e$ (when one identifies $A_e$ with $\lambda_e^\A (A_e)$), see \cite[Propositions 17.13 and 19.3]{Exel2}.

We will need the following lemma several times. It follows readily from \cite[Lemma 17.2]{Exel2}. 
\begin{lemma}\label{E17.2} Assume $a\in (A_e)^+$ and $c \in \A$. Then $c^*ac \in (A_e)^+$ and $c^*ac \leq \|a\| \,c^*c$. Hence,
\[ \|c^*a c\| \, \leq \, \|a\| \, \|c\|^2.\]
\end{lemma}

We will also make use of the matrix C$^*$-algebras associated by Abadie and Ferraro  to a Fell bundle $\A$ in \cite{AF}. We recall here their definition. Let $g_1, \ldots, g_n \in G$ and set $\mathbf{g}:= (g_1, \ldots, g_n) \in G^n$. 
Then
\[ M_\mathbf{g}(\A) := \big\{R=[r_{ij}]\in M_n(\A) : r_{ij} \in A_{g_i^{-1}g_j} \text{ for all } i, j =1, \ldots, n\big\}\]
is a 
$*$-algebra with respect to the natural operations, which can be equipped with a norm turning it into a $C^*$-algebra, 
cf.~\cite[Lemma 2.8]{AF}. 
We note that  $M_\mathbf{g}(\A)$ can also be represented by adjointable operators on a right Hilbert $A_e$-module. Indeed,
consider  the right Hilbert $A_e$-module $A_\mathbf{g}:= A_{g_1^{-1}}\oplus \cdots \oplus A_{g_n^{-1}}$ obtained by taking the direct sum of the $A_{g_i}$'s (considered as right Hilbert $A_e$-modules), whose inner product is given by
\[ \big\langle (a_1, \ldots, a_n), (a'_1, \ldots, a'_n)\big\rangle_{A_e}  = \sum_{i=1}^n a_i^*a'_i\]
for $a_i, a'_i \in A_{g_i^{-1}}$, $i=1, \ldots, n$.
We can let each $R=[r_{ij}] \in M_\mathbf{g}(\A)$ act on $A_\mathbf{g}$ by  
\[L_R(a_1, \ldots, a_n) = (a'_1, \ldots, a'_n), \ \text{ where } a'_i := \sum_{j=1}^n r_{ij}a_j \text { for each }i = 1, \ldots, n. \] 
It is easy to check that the operator $L_R : A_\mathbf{g}\to A_\mathbf{g}$ is adjointable with $(L_R)^*=L_{R^*}$, and that the map $R \mapsto L_R$ from $M_\mathbf{g}(\A)$ into  $\mathcal{L}_{A_e}(A_\mathbf{g})$ is a faithful $*$-homomorphism. 
Thus, the norm on $M_\mathbf{g}(\A)$ satisfies that $\|R\|= \|L_R\|$.

\section{Positive definite bundle maps and complete positivity}
In this section we let $\A=(A_g)_{g\in G}$ be a Fell bundle over a group $G$, $\B=(B_h)_{h\in H}$ be a Fell bundle over a group $H$ and $\varphi:G\to H$ be a group homomorphism, i.e., $\varphi \in {\rm Hom}(G, H)$.

\subsection{Positive definite bundle maps}
\begin{definition}\label{pdbm}
An 
\emph{$\A$-$\varphi$-$\B$ bundle map}
is a family $T = (T_g)_{g\in G}$ of maps $T_g: A_g \to B_{\varphi(g)}$ which are linear and bounded for all $g\in G$. Alternatively, we can think of $T$ as a map $T:\A\to \B$ such that for each $g\in G$ the restriction $T_g$ of $T$ to $A_g$ is a bounded linear map from $A_g$ into $B_{\varphi(g)}$. 

 Note that  $T^*= (T^*_g)_{g\in G}$, defined by
 \[ T^*_g(a) = T_{g^{-1}}(a^*)^* \text{ for all $g\in G$ and $a \in A_g$},\]
 is then also an  $\A$-$\varphi$-$\B$ bundle map.

 An $\A$-$\varphi$-$\B$ bundle map $T$ 
 is said to be 
 \emph{self-adjoint} if it satisfies that  $T^* = T$, i.e.,  \[T_g(a)^* = T_{g^{-1}}(a^*)\] for all $g\in G$ and $a \in A_g$. 
 Moreover, it is said to be \emph{multiplicative} if it satisfies that 
\[ T_{gg'}(aa') = T_g(a) T_{g'}(a')\] for all $g, g' \in G$ and $a\in A_g, a'\in A_{g'}$. 

A pair $(\varphi, T)$ where $\varphi \in {\rm Hom}(G, H)$ and $T$ is  a self-adjoint, multiplicative $\A$-$\varphi$-$\B$ bundle map will be called  a (Fell bundle) {\it morphism from $\A$ to $\B$}.

\end{definition}

 In the case  where $G=H$ and $\varphi = {\rm id}_G$, so $\A$ and $\B$ are both Fell bundles over $G$, we will use the term $\A$-$\B$ \emph{bundle map} instead of $\A$-${\rm id}_G$-$\B$ bundle map. Moreover, if $\A=\B$, we will just say \emph{$\B$-bundle map} instead of $\B$-$\B$ bundle map. 
 
\begin{remark}  Considering the case  where $G=H$ and $\varphi = {\rm id}_G$, Exel defines a morphism from $\A$ to $\B$ to be an $\A$-$\B$ bundle map $T$ which is self-adjoint and multiplicative, and he proves that it gives rise to a $*$-homomorphism from $C^*(\A)$ into $C^*(\B)$ (resp.~from $C_r^*(\A)$ into $C_r^*(\B)$), cf.~\cite[Section 21.1]{Exel2}.

\end{remark}

\begin{definition}
 We will say that an $\A$-$\varphi$-$\B$ 
 bundle map
  $T = (T_g)_{g\in G}$ is \emph{positive definite} whenever we have
\begin{equation}\label{Bposdef}
\sum_{i, j=1}^n b_i \, T_{g_i^{-1}g_j} (a_i^* a_j) \, b_j^* \, \in  (B_e)^+
\end{equation}
for all $n\in \Naturali$, $g_1, \ldots, g_n \in G$ and all  $a_i \in A_{g_i}$, $b_i\in B_{\varphi(g_i)}$, $i=1, \ldots, n$.
\end{definition} 

\begin{remark}
Let  $T = (T_g)_{g\in G}$ be an $\A$-$\varphi$-$\B$ 
 bundle map. We note that $T$ is positive definite if and only if
  for every  $n\in \Naturali$, $\mathbf{g}= (g_1, \ldots, g_n) \in G^n$ and $a_i \in A_{g_i}$, $i=1, \ldots, n$, 
  the matrix 
\[ \Big[ T_{g_i^{-1}g_j} (a_i^* a_j)\Big]\] is positive in the $C^*$-algebra $M_{\varphi(\mathbf{g})}(\B)$, where $\varphi(\mathbf{g}):= (\varphi(g_1), \ldots, \varphi(g_n)) \in H^n$. 

Indeed, consider $\mathbf{g}= (g_1, \ldots, g_n) \in G^n$ and $a_i \in A_{g_i}$ for $i=1, \ldots, n$. 
Let then $S:=[s_{ij}]$ be the matrix in  $M_{\varphi(\mathbf{g})}(\B)$  given by \[s_{ij} := T_{g_i^{-1}g_j} (a_i^* a_j) \in B_{\varphi(g_i^{-1}g_j)} = B_{\varphi(g_i)^{-1}\varphi(g_j)}\] for each $i, j \in \{1, \dots, n\}$. Then for $b_1\in B_{\varphi(g_i)}, \ldots, b_n\in B_{\varphi(g_n)}$, we  have that 
\[ \sum_{i, j=1}^n b_i \, T_{g_i^{-1}g_j} (a_i^* a_j) \, b_j^* = \sum_{i=1}^n b_i\,\Big( \sum_{j=1}^n s_{ij} \, b_j^*\Big) = \big\langle (b_1^*, \ldots, b_n^*), L_S(b_1^*, \ldots, b_n^*) \big\rangle\]
Using \cite[Lemma 4.1]{La1}, it follows from this equality that if $T$ is $\B$-positive definite, then $L_S$ is positive in $\mathcal{L}_{B_e}(B_{\varphi(\mathbf{g})})$. 
As the map $R\mapsto L_R$ is a faithful $*$-homomorphism from $M_{\varphi(\mathbf{g})}(\B)$ into $\mathcal{L}_{B_e}(B_{\varphi(\mathbf{g})})$, this is equivalent to $S$ being positive in $M_{\varphi(\mathbf{g})}(\B)$. This shows the forward implication. The converse implication can be shown by reversing the arguments above.
\end{remark}

\begin{remark} \label{Exel-morphism} 
 Consider  a morphism $(\varphi, T)$ from a Fell bundle $\A=(A_g)_{g\in G}$ to a Fell bundle $\B=(B_h)_{h\in H}$. 
Then $T$ is positive definite. 
Indeed, let $g_1, \ldots, g_n \in G$, $a_i\in A_{g_i}$ and $b_i \in B_{\varphi(g_i)}$ for $i=1, \ldots, n$. Then 
 \begin{align*}
 \sum_{i, j=1}^n b_i \, T_{g_i^{-1}g_j} (a_i^* a_j) \, b_j^* &= \sum_{i, j=1}^n b_i \, T_{g_i^{-1}} (a_i^*) T_{g_j}(a_j) \, b_j^*
  =  \sum_{i, j=1}^n b_i \, T_{g_i} (a_i)^* T_{g_j}(a_j) \, b_j^* 
 \\&=   \Big(\sum_{i=1}^n T_{g_i} (a_i)b_i^*\Big)^* \Big(\sum_{j=1}^n  T_{g_j}(a_j) \, b_j^*\Big) \in  (B_e)^+.
\end{align*}
It will follow from Theorem \ref{cp-full-bundle} that $T$ is full (according to Definition \ref{full-red}). Similarly, if ${\rm ker}(\varphi)$ is amenable, Theorem \ref{cp-reduced} will imply that $T$ is also reduced. The associated maps $\Phi_T: C^*(\A)\to C^*(\B)$ and  $M_T:C_r^*(\A) \to C_r^*(\B)$ are then actually $*$-homomorphisms, as is essentially well-known (cf.~\cite[Propositions 21.2 and 21.3]{Exel2}, where the case $G=H$ and $\varphi = {\rm id}_G$ is considered). 
\end{remark}

\begin{proposition} \label{pd-prop}
Let $T = (T_g)_{g\in G}$ be a positive definite $\A$-$\varphi$-$\B$ bundle map. 
Then the following properties hold.
  \begin{itemize} 
 \item The map $T_e:A_e\to B_e$ is completely positive.
 \item 
 $T$ is self-adjoint.
 \item For $g\in G$ and $a\in A_g$, we have \,$T_g(a)^*T_g(a) \leq \|T_e\| \, T_e(a^*a)$\, $($in $B_e$$)$.
 \end{itemize}
 \end{proposition}
 
  \begin{proof} The first assertion follows readily by choosing all $g_i$'s to be equal to $e$ in the definition of positive definiteness. Next, pick an approximate unit $\{u_i\}_{i\in I} $ for $A_e$. Let $i\in I, g\in G, a\in A_g$. Choosing $g_1 = e, g_2 = g, a_1 = u_i$ and $a_2 = a$, and setting $\mathbf{g}= (e, g)$, we get that  the matrix 
  \[ \begin{bmatrix} T_e(u_i^2) & T_g(u_ia)\\ T_{g^{-1}}(a^*u_i) & T_e(a^*a)\end{bmatrix} \in M_{\varphi(\mathbf{g})}(\B)^+.\]
 In particular, this matrix is self-adjoint, which gives that  $T_g(u_ia)^* = T_{g^{-1}}(a^*u_i)$. 
 
 Taking now the limit w.r.t.~$i$ and using that $T_g$ is bounded, we get that $T_g(a)^* = T_{g^{-1}}(a^*)$, and the second assertion follows. 
 
 Moreover, since $T_e(u_i^2)$ and $T_e(a^*a)$ are positive elements in $A_e$, it follows from \cite[Lemma 5.2 (iii)]{La1}
that\[ T_{g^{-1}}(a^*u_i)T_g(u_ia) = T_g(u_ia)^*T_g(u_ia) \leq \|T_e(u_i^2)\| \,T_e(a^*a) \leq \|T_e(u_i)\| \,T_e(a^*a) .\]
Taking the limit w.r.t.~$i$ and using again that $T_g$ is bounded, we get that the third assertion holds (since $\|T_e\| = \lim_i \|T_e(u_i)\|$ as $T_e$ is completely positive). 
 \end{proof}
 Note that if $\A$ is unital (i.e., if $A_e$ is unital), replacing each $u_i$ with $1_{A_e}$, we do not need to take any limit in the proof above, hence the boundedness assumption for each $T_g$ is not required in this case for the second and the third properties above to hold; the boundedness of each $T_g$ follows then readily from the third property (using that $T_e$ is completely positive, hence bounded), and could therefore be omitted in the definition of a positive definite bundle map.  

As an immediate consequence of Proposition \ref{pd-prop}, we get
\begin{corollary}
Let $T = (T_g)_{g\in G}$ be a positive definite $\A$-$\varphi$-$\B$ bundle map. Then \[\sup_{g\in G} \|T_g\| = \|T_e\|.\]
\end{corollary}

\subsection{Completely positive maps} 

Our interest in positive definite bundle maps is due to their connection with completely positive maps. Here is the first evidence of this fact.

\begin{proposition}\label{cpred-pdmap} Let 
 $M: C_r^*(\A) \to C_r^*(\B)$ be a completely positive linear map. For each $g\in G$, define a  map $T_g:A_g \to B_{\varphi(g)}$ by  
\[T_g(a) = E^\B_{\varphi(g)}\big(M(\lambda^\A_g(a))\big) \quad \text{for each } a\in A_g. \]
Then $T=(T_g)_{g\in G}$ is a positive definite $\A$-$\varphi$-$\B$ bundle map. 
\end{proposition}
\begin{proof} It is clear that each $T_g$ is linear and  bounded, so $T$ is an $\A$-$\varphi$-$\B$ bundle map.  

\noindent Let $g_1, \ldots, g_n \in G$,  $a_i \in A_{g_i}$ and $b_i\in B_{\varphi(g_i)}$ for $i=1, \ldots, n$.
Then for each $i, j$ we have
\[ b_i \, T_{g_i^{-1}g_j} (a_i^* a_j) \, b_j^* = b_i \, E^\B_{\varphi(g_i^{-1}g_j)} \big(M(\lambda^\A_{g_i^{-1}g_j}(a_i^*a_j) \big) \, b_j^* = b_i \, E^\B_{\varphi(g_i)^{-1}\varphi(g_j)} \big( w_{ij} \big) \, b_j^*\]
where $w_{ij} := M(\lambda^\A_{g_i}(a_i)^*\lambda^\A_{g_j}(a_j)) \in C_r^*(\B)$. 

We note that the matrix $W:=[w_{ij}] \in M_n(C_r^*(\B))$ is positive as $M$ is assumed to be completely positive. Considering the direct sum $(\ell^2(\B))^n$
as a right Hilbert $B_e$-module, we let $M_n(C_r^*(\B))$ act on $(\ell^2(\B))^n$ in the natural way. 

Using \cite[Proposition 17.12]{Exel2} and setting 
$\Omega := \big(j^\B_{\varphi(g_1)}(b_1)^*, \ldots, j^\B_{\varphi(g_n)}(b_n)^*\big)  \in (\ell^2(\B))^n$, we therefore get that
\[\sum_{i,j=1}^n b_i \, T_{g_i^{-1}g_j} (a_i^* a_j) \, b_j^* = \sum_{i,j=1}^n  \big\langle j^\B_{\varphi(g_i)}(b_i)^*, w_{ij} \,j^\B_{\varphi(g_j)}(b_j)^*\big\rangle_{\ell^2(\B)} = \big\langle \Omega , W\, \Omega\big\rangle_{(\ell^2(\B))^n}\]
is positive in $B_e$. Thus, $T$ is positive definite. 
\end{proof}
Proceeding in a similar way, it should be straightforward for the reader to check that the following analogous result holds.  

\begin{proposition}\label{cpfull-pdmap}
 Let $\Phi: C^*(\A) \to C^*(\B)$ be a completely positive linear map. Then we get a positive definite $\A$-$\varphi$-$\B$ bundle map $T=(T_g)_{g\in G}$
 by setting 
\[T_g(a) = E^\B_{\varphi(g)}\big(\Lambda^\B(\Phi(\,\widehat{j}^{\,\A}_g(a)))\big) \quad \text{for each } g \in G \text{ and } a\in A_g. \]
\end{proposition}

Let $T = (T_g)_{g\in G}$ be an $\A$-$\varphi$-$\B$ bundle map.  
We can then define a linear map $\phi_T: C_c(\A) \to C_c(\B)$ by setting 
\[  \phi_T\big(f)= \sum_{g\in G} j^\B_{\varphi(g)}(T_g(f(g)))\quad \text{for all } f\in C_c(\A).\]
(Note that the sum above is finite since $f$ has finite support.) It is the only linear map from $C_c(\A)$ to $C_c(\B)$ satisfying that
\[ \phi_T\big(j_g^\A(a)\big)= j^\B_{\varphi(g)}(T_g(a)) \quad \text{for all } g\in G \text{ and } a\in A_g.\]

\begin{definition} \label{full-red}
Let $T = (T_g)_{g\in G}$ be an $\A$-$\varphi$-$\B$ bundle map. We will say that $T$ is \emph{reduced $($resp.~full\,$)$} if $\phi_T$ is bounded w.r.t.~the reduced (resp.~universal) $C^*$-norms on $C_c(\A)$ and $C_c(\B)$, and will denote the natural extension of $\phi_T$ to a linear bounded map from $C_r^*(\A)$ into  $C_r^*(\B)$ (resp.~from $C^*(\A)$ into $C^*(\B)$) by $M_T$ (resp.~$\Phi_T$).
\end{definition}  
   
In other words, $T$ is reduced precisely when there exists a (necessarily unique) linear bounded map  $M_T: C_r^*(\A)\to C_r^*(\B)$ such that 
\[M_T\Big(\lambda^\A_g(a)\Big) = \lambda^\B_{\varphi(g)}(T_g(a)) \, \, \text{ for all} g \in G \text{ and} a \in A_g.\]

Similarly, $T$ is full precisely when there exists a (necessarily unique) linear bounded map $\Phi_T: C^*(\A)\to C^*(\B)$
such that 
 \[\Phi_T\Big(\widehat j^\A_g(a)\Big) =  \widehat j^\B_{\varphi(g)}(T_g(a)) \, \, \text{ for all} g \in G \text{ and} a \in A_g.\]

We note that if $T$ is reduced (resp. full) and $M_T$ (resp. $\Phi_T$) is completely positive, then $T$ is positive definite. 
Indeed, the first claim follows from Proposition \ref{cpred-pdmap} (with $M = M_T$), since  
\[T_g(a) = E^\B_{\varphi(g)}\big(\lambda^\B_{\varphi(g)}(T_g(a))\big) = E^\B_{\varphi(g)}\big(M_T(\lambda^\A_g(a))\big) \, \, \text{ for all} g \in G \text{ and} a \in A_g,\] 
while the second one follows in a similar fashion from Proposition \ref{cpfull-pdmap} (with $\Phi= \Phi_T$).

\medskip
We will now discuss whether the converse of these last two assertions hold. We will first 
show that it does in the full case. 
To obtain a reduced version, we will have to assume that the kernel of $\varphi$ is amenable. Such a condition 
is known to be necessary and sufficient when considering the extendibility of the homomorphism $\varphi:G\to H$ to a $*$-homomorphism from $C_r^*(G)$ into $C_r^*(H)$
 (cf.~\cite[Corollary 8.C.15]{BedlH}). Note here that a completely positive map from $C_r^*(G)$ into $C_r^*(H)$ which extends $\varphi$ has to be a $*$-homomorphism, so one gets from the result cited above that $\varphi$ extends to a completely 
positive linear map from $C_r^*(G)$ into $C_r^*(H)$  if and only if ${\rm ker}(\varphi)$ is amenable.

\begin{theorem}  \label{cp-full-bundle} Let $T=(T_g)_{g\in G}$ be an $\A$-$\varphi$-$\B$ bundle map. Then $T$ is positive definite if and only if  $T$ is full and the associated map $\Phi_T:C^*(\A)\to C^*(\B)$ is completely positive. In such a case, we have $\|\Phi_T\| = \|T_e\|$.
\end{theorem}

\begin{proof} 
In view of our preliminary discussion, in order to prove the first statement, we only need to show the forward implication.
 So assume that $T$ is positive definite. 
Let $\psi_T : \A \to C^*(\B)$ be the map given by \[\psi_T(a) = \widehat j^\B_{\varphi(g)}(T_g(a))\]
 whenever $a\in A_g$ for $g\in G$,
and let $\Pi$ be a faithful $*$-representation of $C^*(\B)$ on some Hilbert space $\H$. 

As a first step, we will show that $\phi:=\Pi\circ \psi_T: \A \to B(\H)$ is completely positive in the sense of Buss, Ferraro and Sehnem, cf.~\cite[Appendix]{BFS}.  
Let $g_1, \ldots, g_n \in G$, set $\mathbf{g}:= (g_1, \ldots, g_n) \in G^n$ 
and let $\phi^\mathbf{g} : M_\mathbf{g}(\A) \to M_n(B(\H))$ be the linear map defined by 
\[ \phi^\mathbf{g}([r_{ij}]) = [\phi(r_{ij})] = \big[\Pi(\psi_T(r_{ij}))\big] \, \text{ for each } [r_{ij}] \in M_\mathbf{g}(\A).\]
To show that $\phi$ is completely positive amounts to check that $\phi^\mathbf{g}$ is a positive map for every $\mathbf{g} \in \cup_{n\in \Naturali} \,G^n$. It clearly suffices to show that 
the linear map $\psi^\mathbf{g}_T: M_\mathbf{g}(\A)\to M_n(C^*(\B))$, given by \[\psi^\mathbf{g}_T([r_{ij}]) = \big[\psi_T(r_{ij})\big] = \Big[ \,\widehat j^\B_{\varphi(g_i^{-1}g_j)}(T_{g_i^{-1}g_j}(r_{ij}))\Big] \, \text{ for each } R= [r_{ij}] \in M_\mathbf{g}(\A), \]
is a positive map for every $\mathbf{g} \in \cup_{n\in \Naturali} \,G^n$. 
So let $\mathbf{g}= (g_1, \ldots, g_n) \in G^n$ and consider $a_i \in A_{g_i}$ for $i=1, \ldots, n$. Since $T$ is positive definite, the
matrix $\big[T_{g_i^{-1}g_j}(a_{i}^*a_j)\big]$ is positive in $M_{\varphi(\mathbf{g})}(\B)$.
 Now, it is easy to check that if $\mathbf{h}= (h_1, \ldots, h_n) \in H^n$, then the map 
 \[\big[s_{ij}\big] \in  M_\mathbf{h}(\B) \mapsto \Big[\, \widehat{j}^\B_{h_i^{-1}h_j} (s_{ij})\Big] \in M_n(C^*(\B))\] is a $*$-homomorphism. Hence we get 
 that 
 \[ \psi^\mathbf{g}_T\big([a_i^*a_j]\big) = \Big[ \,\widehat j^\B_{\varphi(g_i)^{-1}\varphi(g_j)}(T_{g_i^{-1}g_j}(a_{i}^*a_j))\Big] \geq 0\,.\]
 It then follows 
 that $\psi_T^\mathbf{g}([r_{ij}]) \geq 0$ for every $R=[r_{ij}] \in M_\mathbf{g}(\A)^+$, hence that $\psi_T^\mathbf{g}$ is a positive map. 
 
 Indeed, assume $R= C^*C$ for some $C= [c_{ij}] \in M_\mathbf{g}(\A)$. Then 
 $R  = \sum_{k=1}^n \big[ c_{ki}^*c_{kj}\big],$
 so applying what we have just shown gives that \[\psi^\mathbf{g}_T(R) = \sum_{k=1}^n \psi^\mathbf{g}_T\big(\big[ c_{ki}^*c_{kj}\big]\big)\] is positive in $M_n(C^*(\B))$, being a finite sum of positive elements.
 
 This shows that $\phi=\Pi\circ \psi_T: \A \to B(\H)$  is completely positive, as desired.

\medskip \noindent By \cite[Corollary A.8]{BFS}, there exists a completely positive map $\tilde\phi: C^*(\A) \to B(\H)$ (called the ``integrated form'' of $\phi$ 
in \cite{BFS}) such that
\[ \tilde\phi\big(\widehat j^\A_g(a)\big) =\phi(a)\quad \text{for every } g \in G \text{ and } a \in A_g.\] 
We then get that 
\[   \tilde\phi(\widehat j^\A_g(a)) =\phi(a) = \Pi(\psi_T(a)) = \Pi\big(\,\widehat j^\B_{\varphi(g)}(T_g(a))\big) \quad \text{for every } g \in G \text{ and } a \in A_g.\]
Thus, $\tilde\phi\big(C^*(\A)\big) \subseteq \Pi(C^*(\B))$. Since $\Pi$ is injective, we may regard it as a $*$-isomorphism onto its range $\Pi(C^*(\B))$ and define $\phi' :C^*(\A) \to C^*(\B)$ by 
$\phi' = \Pi^{-1}\circ \tilde\phi$. Then $\phi'$ is completely positive and satisfies
\[ \phi'\big(\widehat j^\A_g(a)\big) = \widehat j^\B_{\varphi(g)}(T_g(a))\quad \text{for every } g \in G \text{ and } a \in A_g.\] 
Thus, $\phi'$ extends $\phi_T$, so $T$ is full, and $\Phi_T = \phi'$ is completely positive, as we wanted to show.

Finally,  let $\{u_i\}_{i \in I}$ be an approximate unit for $A_e$. Then $\big\{\widehat j_e^\A(u_i)\big\}_{i\in I}$ is an approximate unit for $C^*(\A)$. So if $T$ is positive definite, then $T_e$ and $\Phi_T$ are completely positive, and  thus
\[ \|\Phi_T\| = \lim_i \|\Phi_T(\,\widehat j_e^\A(u_i))\|_{\rm u} =  \lim_i \|\widehat j_e^\A(T_e(u_i))\|_{\rm u} =  \lim_i \|T_e(u_i)\| = \|T_e\|.\]
 
 \vspace{-2ex}\end{proof}
Note that if $\A$ and $\B$ are amenable in the sense of Exel \cite{Exel1, Exel2}, i.e., their full and reduced cross-sectional $C^*$-algebras are  canonically $*$-isomorphic, respectively, then we get as a corollary to Theorem \ref{cp-full-bundle} that $T$ is a positive definite $\A$-$\varphi$-$\B$ bundle map if and only if  $T$ is reduced and $M_T:C_r^*(\A)\to C_r^*(\B)$ is completely positive. 
 When $\A$ or $\B$ are not amenable, this equivalence is also true under the assumption that ${\rm ker}(\varphi)$ is an amenable subgroup of $G$. In order to show this, we will first generalize \cite[Proposition 4.8]{BFS}. We will use the following terminology.
 \begin{definition}
 A linear map $\Phi:C^*(\A)\to C^*(\B)$ (resp.~$\Psi:C_r^*(\A)\to C_r^*(\B)$) will be called \emph{$\varphi$-diagonal} whenever
 $\Phi(\,\widehat{j}^{\,\A}_g(A_g)) \subseteq \,\,\widehat{j}^{\,\B}_{\varphi(g)}(B_{\varphi(g)})$
 (resp.~$\Psi\big(\lambda^\A_g(A_g)\big) \subseteq \,\lambda^\B_{\varphi(g)}(B_{\varphi(g)})$)
for all $g\in G$. 
 \end{definition}
 For example, if $T$ is a full (resp.~reduced) $\A$-$\varphi$-$\B$ bundle map, then the associated map $\Phi_T$ (resp.~$M_T$) is $\varphi$-diagonal. 
 It is easy to see that  $\{\Phi_T: T \text{ is a full $\A$-$\varphi$-$\B$ bundle map}\}$ coincides with the set of $\varphi$-diagonal bounded linear maps from $C^*(\A)$ into $C^*(\B)$.  A similar statement holds for $\varphi$-diagonal bounded linear maps from $C_r^*(\A)$ into $C_r^*(\B)$. 

 \begin{proposition}\label{BFS_mod} Assume that the kernel of $\varphi:G\to H$ is amenable  and let $\pi = (\pi_h)_{h\in H}$ be a 
 representation of $\B$ on a Hilbert space $\H$ such that $\pi_e$ is faithful. Suppose that $\psi: C^*(\A) \to B(\H)$ is a completely positive linear map such that  for all $g\in G$ \[\psi(\,\widehat j^\A_g(A_g)) \subseteq \, \pi_{\varphi(g)}(B_{\varphi(g)}).\] 
  Then there is a $\varphi$-diagonal completely positive linear map $\psi': C^*_r(\A) \to C_r^*(\B)$ 
 such that  for all $g\in G$
 \begin{equation}\label{psiprime2}
 \big(\pi_{\varphi(g)}\circ(\lambda^\B_{\varphi(g)})^{-1}\circ\psi'\circ{\lambda^\A_g}\big)
 _{|A_g} 
 = \big(\psi\circ\widehat j^\A_g\big)
 _{|A_g}.
 \end{equation}
  \end{proposition} 
\begin{proof}
By Stinespring's dilation theorem for completely positive maps (see e.g.~\cite[Theorem 5.6]{La1}), there exist a Hilbert space $\K$, a $*$-homomorphism $\widehat \psi: C^*(\A) \to B(\K)$ and a bounded operator $V: \H \to \K$ such that
$\psi(x) = V^*\,\widehat\psi(x)\,V$
for all $x \in C^*(\A)$. Let $\widehat\pi=(\widehat\pi_g)_{g\in G}$ denote the representation of $\A$ on $\K$ associated to $\widehat \psi$, so that 
\[ \widehat\psi (\,\widehat j_g^\A(a))= \widehat\pi_g(a)\quad \text{ for all } g \in G \text{ and } a\in A_g.\]

By Fell's absorption principle for Fell bundles \cite[Proposition 18.4]{Exel2}, the integrated forms $\phi_\A$ and $\phi_\B$ of 
the representations 
 $(\widehat\pi_g\otimes\lambda^G_g)_{g\in G}$ and 
 $(\pi_h\otimes\lambda^H_h)_{h\in H}$ (as defined in \cite[Prop~16.29]{Exel2}) factor through $C^*_r(\A)$ and $C^*_r(\B)$, respectively. This means that there are $*$-homomorphisms  
 $\psi_\A : C^*_r(\A) \to B(\K\otimes \ell^2(G))$ and $\psi_\B : C^*_r(\B) \to B(\H \otimes \ell^2(H))$ 
 such that $\phi_\A=  \psi_\A\circ\Lambda^\A$ and $\phi_\B=  \psi_\B\circ\Lambda^\B$. 
 Using \cite[Corollary 18.5]{Exel2},  we may regard $\psi_\A$ as a map from $C^*_r(\A)$ into $B(\K)\otimes_{\rm min}C_r^*(G)$ and $\psi_\B$ as a map from $C^*_r(\B)$ into $B(\H)\otimes_{\rm min}C_r^*(H)$, which satisfy that
\[ \psi_\A(\lambda_g^\A(a)) = \widehat{\pi}_g(a)\otimes\lambda_g^G \, \text{ and } \, \psi_\B(\lambda_h^\B(b)) = \pi_h(b)\otimes\lambda_h^H\]
for all $g\in G$, $a\in A_g$, $h \in H$ and $b\in B_h$. 
Note that  $\psi_\B$ is faithful since $\pi_e$ is assumed to be so, cf.~\cite[Proposition 18.4]{Exel2}. 

Now, as shown in \cite[Corollary 8.C.15]{BedlH}, since $\ker(\varphi)$ is amenable, $\varphi$ extends to a $*$-homomorphism $\varphi_{*, r}: C_r^*(G)\to C_r^*(H)$ satisfying that $\varphi_{*, r}(\lambda_g^G) = \lambda^H_{\varphi(g)}$ for all $g\in G$.

Let $\widetilde V: B(\K)\to B(\H)$ denote the completely positive map given by $\widetilde V(T) = V^*TV$ for all $T\in B(\K)$. 
Using \cite[Theorem 3.5.3]{BrOz}, we may define a completely positive map $\Phi:  C_r^*(\A) \to B(\H)\otimes_{\rm min}C_r^*(H)$ by 
\[ \Phi:= (\widetilde V \otimes_{\rm min} \varphi_{*, r})\circ \psi_\A.\]
Let $g \in G$ and $a \in A_g$. Then we have
\begin{align*} 
\Phi(\lambda^\A_g(a)) & = (\widetilde V \otimes_{\rm min} \varphi_{*, r})(\psi_\A(\lambda^\A_g(a)))   =  (\widetilde V \otimes_{\rm min} \varphi_{*, r})(\widehat\pi_g(a) \otimes \lambda^G_g)  \\
& = V^*\widehat\pi_g(a)V \otimes  \lambda^H_{\varphi(g)} = V^* \widehat\psi(\widehat j^\A _g(a)) V \otimes \lambda^H_{\varphi(g)} \\
&= \psi(\,\widehat j^\A _g(a))  \otimes \lambda^H_{\varphi(g)}.
\end{align*}
By assumption, we have that $\psi(\,\widehat j^\A_g(a)) \in  \pi_{\varphi(g)}(B_{\varphi(g)})$. Thus there is some $b \in B_{\varphi(g)}$ such that $\psi(\,\widehat j^\A_g(a)) = \pi_{\varphi(g)}(b)$. Hence 
\[\Phi(\lambda^\A_g(a)) = \pi_{\varphi(g)}(b)\otimes 
\lambda^H_{\varphi(g)} = \psi_\B(\lambda^\B_{\varphi(g)}(b)).\]
By a density and continuity argument, we obtain that $\Phi(C_r^*(\A)) \subseteq \psi_\B(C_r^*(\B))$.

Since $\psi_\B$ is faithful, we can define  
$\psi': C_r^*(\A)\to C_r^*(\B)$  by $\psi' = (\psi_\B)^{-1}\circ \Phi$, which then satisfies 
$\psi'(\lambda^\A_g(A_g))\subseteq \lambda^\B_{\varphi(g)}(B_{\varphi(g)})$, i.e., $\psi'$ is $\varphi$-diagonal.
It is now straightforward to check that $\psi'$ also satisfies (\ref{psiprime2}).
\end{proof}

\begin{theorem}  \label{cp-reduced} Let $T=(T_g)_{g\in G}$ be an $\A$-$\varphi$-$\B$ bundle map, and assume that the kernel of $\varphi$ is amenable. Then $T$ is positive definite if and only if  $T$ is reduced  and the associated map $M_T:C_r^*(\A)\to C_r^*(\B)$ is completely positive, in which case we have
$\|M_T\| 
= \|T_e\|$.
\end{theorem}
\begin{proof} As we have seen previously, the backward implication is always true. So assume that $T=(T_g)_{g\in G}$ is positive definite. By Theorem \ref{cp-full-bundle} we know that $\phi_T$ extends to a completely positive map $\Phi_T: C^*(\A)\to C^*(\B)$ satisfying that $\Phi_T\big(\widehat j^\A_g(a)\big) = \widehat j^\B_{\varphi(g)}(T_g(a))$ for all $g\in G$ and $a\in A_g$. 
Let $\pi=(\pi_h)_{h\in H}$ be a $*$-representation of $\B$ in $B(\H)$ for some Hilbert space $\H$ such that 
its integrated form $\Pi:C^*(\B) \to B(\H)$ 
is faithful. Note that $\pi_e= \Pi \circ \widehat j^{\,\B}_e$ is faithful too. 
Then the linear map $\widetilde\phi :=\Pi \circ \Phi_T:C^*(\A) \to B(\H)$ is completely positive and satisfies that 
\[\widetilde\phi\Big(\widehat j^\A_g(a)\Big) = \Pi\big(\Phi_T(\widehat j^\A_{g}(a))\big) 
= \Pi\big(\widehat j^\B_{\varphi(g)}(T_g(a))\big)  = 
\pi_{\varphi(g)}(T_g(a))\]
for all $g\in G$ and $a\in A_g$. Now, using Proposition \ref{BFS_mod},
we get that there exists a completely positive map $\phi':C_r^*(\A)\to C_r^*(\B)$ such that $\phi'\big(\lambda^\A_g(A_g)\big) \subseteq \lambda^\B_{\varphi(g)}(B_{\varphi(g)})$ and 
\[\big(\pi_{\varphi(g)}\circ(\lambda^\B_{\varphi(g)})^{-1}\circ\phi'\circ{\lambda^\A_g}\big)_{|A_g} = \big(\widetilde\phi\circ\widehat j^\A_g\big)_{|A_g}\] for each $g\in G$. This means that
\[ \pi_{\varphi(g)}\big((\lambda^\B_{\varphi(g)})^{-1}( \phi'(\lambda^\A_g(a)))\big)  = \pi_{\varphi(g)}(T_g(a))\]
for $g\in G$ and $a\in A_g$. Since $\pi_{\varphi(g)}$ is injective on $B_{\varphi(g)}$ 
(this readily follows from the definition of a $*$-representation \cite[Def. 16.20]{Exel2} and the fact that $\pi_e$ is faithful by assumption), this implies that 
\[ \phi'(\lambda^\A_g(a))=\lambda^\B_{\varphi(g)}(T_g(a))\quad\text{for every } g\in G \text{ and } a\in A_g.\]
This shows that $\phi':C_r^*(\A)\to C_r^*(\B)$ is a bounded extension of $\phi_T$. Hence $T$ is reduced, and $M_T = \phi'$ is completely positive, as desired. 

Finally, if $T$ is positive definite, then one may proceed in a similar way as in the full case to show that $\|M_T\| = \|T_e\|$.
\end{proof}

\begin{remark} 
Even if ${\rm ker}(\varphi)$ is not amenable it may still happen that a positive definite $\A$-$\varphi$-$\B$ bundle map $T$ is reduced with $M_T$ completely positive, cf.~Example \ref{group bundle}. 
Hence, in such a setting, 
 it would be interesting to find some additional condition(s)  ensuring that $T$ is reduced with $M_T$ completely positive. 
\end{remark}

\begin{remark}
Given a completely positive linear map $\Phi:C^*(\A)\to C^*(\B)$, let $T^\Phi$ denote the positive definite $\A$-$\varphi$-$\B$ bundle map associated to it in Proposition \ref{cpfull-pdmap}.  Then, by Theorem \ref{cp-full-bundle}, the map $\Phi_{T^\Phi}$ is a $\varphi$-diagonal completely positive linear map (which could be called the \emph{$\varphi$-diagonalization of $\Phi$}). Hence the map $\Phi\mapsto  \Phi_{T^\Phi}$ is a projection from the space CP($C^*(\A), C^*(\B)$) of completely positive linear maps onto the subspace of
$\varphi$-diagonal maps.  
If the kernel of $\varphi$ is amenable, using Proposition \ref{cpred-pdmap} and Theorem \ref{cp-reduced},  one obtains a similar  projection map on $CP(C_r^*(\A),C_r^*(\B))$.
\end{remark}

\begin{remark} \label{FS}
The family of $\A$-$\varphi$-$\B$ bundle maps forms a vector space $L(\A, \varphi, \B)$ with respect to pointwise addition and scalar multiplication. Set 
\[B(\A, \varphi, \B):={\rm span}\{T \in L(\A, \varphi, \B): T \text{ is positive definite}\}.\] Then any $S \in B(\A, \varphi, \B)$ is a full $\A$-$\varphi$-$\B$ bundle map, which is also reduced if ${\rm ker}(\varphi)$ is assumed to be amenable. This follows readily from Theorems \ref{cp-full-bundle} and \ref{cp-reduced} after noticing that the map $T \mapsto \phi_T$ is linear. We thus get a linear map $T\mapsto \Phi_T$ from $B(\A, \varphi, \B)$ into the space of completely bounded linear maps CB($C^*(\A), C^*(\B)$),  and a similar one into CB($C_r^*(\A), C_r^*(\B)$) when ${\rm ker}(\varphi)$ is amenable. It is easy to see that
the space $B(\A) := B(\A, {\rm id}_G, \A)$ becomes a unital algebra w.r.t.~natural composition, which could be called the \emph{Fourier-Stieltjes algebra of $\A$}, in analogy with 
 the case where $\A$ is the Fell bundle associated to a unital discrete twisted $C^*$-dynamical system, cf.~\cite{BeCo6}.  
\end{remark}   

When $G=H$ and $\varphi = {\rm id}_G$, we may combine Theorem \ref{cp-full-bundle} and Theorem \ref{cp-reduced} to obtain the following. 
\begin{corollary} \label{cp-full-reduced-cor} Assume $\A$ and $\B$ are both Fell bundles over $G$ and let $T=(T_g)_{g\in G}$ be an $\A$-$\B$ bundle map. Then we have that $T$ is positive definite if and only if $T$ is full  and the associated map $\Phi_T:C^*(\A)\to C^*(\B)$ is completely positive, if and only if $T$ is reduced  and the associated map $M_T:C_r^*(\A)\to C_r^*(\B)$ is completely positive. When this happens, we have
$\|\Phi_T\| = \|M_T\| 
= \|T_e\|$. 
\end{corollary}

\begin{remark}\label{category}
 One can define a category ${\tt Fell}$-${\tt PD}$ whose objects  are Fell bundles over discrete groups, and 
arrows from $\A=(A_g)_{g\in G}$ to $\B=(B_h)_{h\in H}$ are pairs $(\varphi, T)$, where
$\varphi:G\to H$ is a homomorphism and $T$ is a positive definite $\A$-$\varphi$-$\B$ bundle map, composition of arrows being defined in the natural way. 

We then get a functor
$\mathcal{F}$ from the category ${\tt Fell}$-${\tt PD}$ to the category ${\tt C^*CP}$ whose objects are $C^*$-algebras and  arrows are completely positive linear maps,  by setting 
\[ \mathcal{F}(\A) := C^*(\A) \, \text{ and } \, \mathcal{F}(\varphi, T) := \Phi_T:C^*(\A)\to C^*(\B) \]
for an object $\A$ in ${\tt Fell}$-${\tt PD}$ and an arrow $(\varphi, T)$ from $\A$ to  $\B$ in ${\tt Fell}$-${\tt PD}$.

Moreover, ${\tt Fell}$-${\tt PD}$ contains a subcategory ${\tt Fell}$ with the same objects, but where 
arrows are
required to be morphisms according to Definition \ref{pdbm}. 
As $\Phi_T = \mathcal{F}(\varphi, T)$ 
is then a $*$-homomorphism from $C^*(\A)$ into $C^*(\B)$,
we may restrict $\mathcal{F}$ to ${\tt Fell}$ and  get a functor 
into the category ${\tt C^*HOM}$ of $C^*$-algebras having $*$-homomorphisms as arrows.

Similarly, we may also consider the subcategory ${\tt Fell}$-${\tt PD}_{\tt red}$ of ${\tt Fell}$-${\tt PD}$ having the same objects but whose arrows involve only group homomorphisms with amenable kernels,  
and obtain a functor $\mathcal{F}_{\rm red}$ from ${\tt Fell}$-${\tt PD}_{\tt red}$ to ${\tt C^*CP}$.
By restricting $\mathcal{F}_{\rm red}$ to 
the subcategory ${\tt Fell}_{\tt red}$  having the same objects, but whose arrows are morphisms $(\varphi, T): \A\to \B$ 
with ${\rm  ker}(\varphi)$ is amenable,
we  get a functor from ${\tt Fell}_{\tt red}$ into  ${\tt C^*HOM}$.

\end{remark}

\section{Some examples}
\begin{example}\label{group bundle}
We first consider the simple case where $\A = (\Complessi \times \{g\})_{g\in G}$ is the group bundle associated to $G$, $\B= \Complessi \times \{e\}$ is the trivial bundle over $\{e\}$
and $\varphi_0: G\to \{e\}$ is the trivial homomorphism. Then one readily sees that  an $\A$-$\varphi_0$-$\B$ 
 bundle map $T=(T_g)_{g\in G}$ is of the form $T_g(z,g) = (\phi(g)z, e)$ for some (necessarily bounded) function $\phi:G\to \Complessi$, and that $T$ is positive definite if and only if $\phi$ is positive definite.  In this case, we have $C^*(\A) = C^*(G)$, $C^*(\B) = \Complessi$, and the map $T\mapsto \Phi_T$ corresponds to the usual correspondence between positive definite functions on $G$ and positive linear functionals on $C^*(G)$. 
 
 Moreover, we have $C_r^*(\A) = C_r^*(G)$, $C_r^*(\B) = \Complessi$. 
 Thinking of a positive definite $\A$-$\varphi_0$-$\B$ 
 bundle map $T$ as a positive definite function $\phi$ on $G$, 
 it follows for example from \cite[Lemma 2.1.2]{KL} that $M_T$ exists as a positive linear functional on $C_r^*(G)$  exactly when the unitary representation of $G$ associated with $\phi$  via the GNS (or Gelfand-Raikov) construction is weakly contained in $\lambda^G$ (in the sense of \cite{Dix, KL}). 
  Note that this is true regardless of whether $G= {\rm ker} (\varphi_0)$ is amenable or not.
 
 \end{example}

\begin{example}
 Next, if both $\A$ and $\B$ are Fell bundles over the trivial group $\{e\}$ with fibers $A$ and $B$,
  respectively, 
  then a positive definite $\A$-$\B$ bundle map $T$ corresponds 
 to a completely positive map from $A$ into $B$. Of course, in this case, we have $C^*(\A)=A$ and $C^*(\B)=B$. 
\end{example}

\begin{example} 
 We now look at a more general situation, covering both previous examples. Let $\A= (A_g)_{g\in G}$ be a Fell bundle and $B$ be a $C^*$-algebra. Let $\B$ be the Fell bundle over the trivial group $\{e\}$ with fiber $B$ over $e$, so $B= C^*(\B)$, and  let $\varphi_0:G \to \{e\}$ be the trivial homomorphism. Then a positive definite $\A$-$\varphi_0$-$\B$ bundle map $T=(T_g)_{g\in G}$ is a family of linear maps $T_g: A_g \to B$ such that the matrix $\big[ T_{g_i^{-1}g_j} (a_i^* a_j)\big]$ is  positive in $M_n(B)$ for every $g_1, \ldots, g_n \in G$ and $a_i \in A_{g_i}, i=1, \ldots, n$. 
 
In particular, choosing $B=\Complessi$, we get that a positive definite $\A$-$\varphi_0$-$\Complessi$ bundle map $\omega=(\omega_g)_{g\in G}$ is a family of bounded linear functionals $\omega_g: A_g \to \Complessi$ such that the matrix $\big[ \omega_{g_i^{-1}g_j} (a_i^* a_j)\big]$ is  positive in $M_n(\Complessi)$ for every $g_1, \ldots, g_n \in G$ and $a_i \in A_{g_i}, i=1, \ldots, n$. This is easily seen to be equivalent to requiring that
\[ \sum_{i, j=1}^n \omega_{g_i^{-1}g_j} (a_i^* a_j) \, \geq \, 0\]
for every $g_1, \ldots, g_n \in G$ and $a_i \in A_{g_i}, i=1, \ldots, n$. 
Such a bundle map $\omega=(\omega_g)_{g\in G}$ is essentially the same as what Fell and Doran call a 
 functional on $\A$ of positive type, cf.~\cite[Chapter VIII, 20.2]{FD2}. We refer to \cite{AG} for a recent study of such maps. In accordance with Theorem \ref{cp-full-bundle},
every such a functional $\omega=(\omega_g)_{g\in G}$ on $\A$ gives rise to a positive linear functional $\Phi_\omega$ on $C^*(\A)$  satisfying that 
 $\Phi_\omega(\,\widehat j_g^\A(a)) = \omega_g(a)$  for $g \in G$ and $a\in A_g$. Conversely, each 
positive linear functional $\Psi$ on $C^*(\A)$ gives rise to a positive definite $\A$-$\varphi_0$-$\Complessi$ bundle map $\omega^\Psi=(\omega^\Psi_g)_{g\in G}$ given by $\omega_g^\Psi(a) = \Psi(\,\widehat j_g^\A(a))$ for $g \in G$ and $a\in A_g$, which satisfies that
$\Psi = \Phi_{\omega^\Psi}$. It follows that 
the map $\omega \mapsto \Phi_\omega$ is a bijection between the cone of positive definite $\A$-$\varphi_0$-$\Complessi$ bundle maps and the cone of positive linear functionals on $C^*(\A)$.

In another direction, let now $\A = (\Complessi \times \{g\})_{g\in G}$ be the group bundle over $G$ and  assume $B =B(\H)$ for some Hilbert space $\H$.  Then 
a positive definite $\A$-$\varphi_0$-$\B$ bundle map amounts to a map $T: G\to B(\H)$ which is completely positive definite in the sense  of Paulsen \cite[Chapter 4, p.~51]{Pau} (some other authors just say positive definite, e.g.~\cite[Section 6, p.~115]{Bo} or \cite[Section 8.3]{vNe}), i.e., it satisfies that the matrix $[T(g_i^{-1}g_j)]$ is positive in $M_n(B(\H))$ for every $g_1, \ldots, g_n \in G$.  As mentioned by Paulsen on p.~54, there is a one-to-one correspondence between such completely positive definite maps from $G$ into $B(\H)$ and completely positive maps from $C^*(G)$ into $B(\H)$. This fits with Theorem \ref{cp-full-bundle}. 
\end{example}

\begin{example} \label{pd-function} Assume $\B=(B_g)_{g\in G}$ is a Fell bundle.
Let $\phi: G\to \Complessi$ be positive definite
and let $T^\phi= (T^\phi_g)_{g\in G}$ be the $\B$-bundle map given by 
\[ T^\phi_g(b) = \phi(g) b\]
for every $g\in G$ and $b\in B_g$. 
Then $T^\phi$ is positive definite. Indeed, for $g_1, \ldots, g_n \in G$ and $b_i, c_i \in B_{g_i}$, $i=1, \ldots, n$, we have
\[ \sum_{i, j=1}^n c_i \, T^\phi_{g_i^{-1}g_j}(b_i^* b_j) \, c_j^* =  \sum_{i, j=1}^n c_i \, \phi(g_i^{-1}g_j)b_i^* b_j \, c_j^*
=  \sum_{i, j=1}^n  \phi(g_i^{-1}g_j)\, a_{ij}\]
where $a_{ij}:=  c_i b_i^* b_j  c_j^* \in B_e$. Now, the matrix $[ \phi(g_i^{-1}g_j)]$ is positive in $M_n(\Complessi)$ since $\phi$ is positive definite, and $[a_{ij}]$ is clearly positive in $M_n(B_e)$. Hence, the matrix  $[ \phi(g_i^{-1}g_j)a_{ij}]$ is then positive in $M_n(B_e)$ (cf.~\cite[Lemma 3.1]{DoRu}) and the claim readily follows. 

Thus, Theorems \ref{cp-full-bundle} and \ref{cp-reduced} apply and give completely positive linear maps $\Phi^\phi: C^*(\B)\to C^*(\B)$ and $M^\phi: C_r^*(\B)\to C_r^*(\B)$ determined by 
\[ \Phi^{\phi}(j_g^{\B}(b)) = \phi(g) j_g^{\B}(b), \quad M^{\phi}(\lambda_g^{\B}(b)) = \phi(g) \lambda_g^{\B}(b)\]
for all $g\in G$ and $b\in B_g$.

Note that if $\gamma:G\to \Toro$ is a character, then $T^\gamma$ is easily seen to be a morphism from $\B$ into itself, and the associated  map $\Phi^\gamma$ (resp.~$M^\gamma$) is a $*$-automorphism of $C^*(\B)$ (resp.~$C_r^*(\B)$).   
\end{example}

\begin{example} \label{Exel-mult} 
Let  $\B=(B_g)_{g\in G}$ be Fell bundle.
We may then consider maps of the type used by Exel to introduce the approximation property for Fell bundles \cite{Exel1, Exel2}, see also Example \ref{AP-BC}.
Let $\xi:G\to B_e$ be a finitely supported function. For each $g\in G$, define a linear map $T_g: B_g\to B_g$ by
 \[ T_g (b) = \sum_{h\in G} \xi(gh)^* b \, \xi(h) \quad \text{for all } b\in B_g.\]
Then $T = (T_g)_{g\in G}$ is a positive definite $\B$-bundle map.
Indeed, let $g_1, \ldots, g_n \in G$ and  $b_i, c_i \in B_{g_i}$ for $i=1, \ldots, n$. For each $j \in \{1, \ldots , n\}$, let $d_j: G \to B_e$ be the finitely supported function defined by \[d_j(g) = b_j \,\xi(g_j^{-1}g)\, c_j^*\] for all $g\in G$, and consider  $d_j$ as a function in the Hilbert $B_e$-module $\ell^2(G, B_e)$. Then for each $i, j \in \{1, \ldots , n\}$ we have 
\[ c_i \, T_{g_i^{-1}g_j} (b_i^* b_j) \, c_j^*= \sum_{h \in G} c_i\, \xi(g_i^{-1}g_j h)^* b_i^*b_j\, \xi(h) \,c_j^*= \sum_{g\in G} d_i(g)^*d_j(g) = \langle d_i, d_j\rangle_{B_e}.\]
Hence, 
\[ \sum_{i, j=1}^n c_i \, T_{g_i^{-1}g_j} (b_i^* b_j) \, c_j^* = \sum_{i, j=1}^n\langle d_i, d_j\rangle_{B_e}
 = \big\langle \sum_{i=1}^n d_i, \, \sum_{j=1}^n d_j\big\rangle_{B_e} \in (B_e)^+. \]
\end{example}

\begin{example}  \label{dynsyst} 
Let $\A =(A_g)_{g\in G}$ and $\B=(B_h)_{h\in H}$ denote the canonical Fell bundles associated to some discrete $C^*$-dynamical systems $\Sigma= (A,  G, \alpha)$ and $\Omega=(B, H, \beta)$. For simplicity, we only consider untwisted systems here, although the case where both $\Sigma$ and $\Omega$ are twisted systems 
can be handled in a similar way.  We recall that each fibre $A_g := \{ (a, g): a\in A\} $ is identified with $A$ as a Banach space via the map $(a, g) \mapsto a$, the multiplication $A_g\times A_{g'} \to A_{gg'}$ and the involution $A_g\to A_{g^{-1}}$ being given by \[(a, g) (a', g') := (a\alpha_g(a'), gg'), \text{ and } 
 (a,g)^* := (\alpha_{g^{-1}}(a^*), g^{-1}),\] 
 respectively. Each fiber  $B_h$ of $\B$ 
 is defined similarly, along with the Fell bundle operations. 
 
 We also recall that 
 $C^*(\A)$ (resp.~$C_r^*(\A)$) corresponds to the full (resp.~reduced) $C^*$-crossed product $C^*(\Sigma)$ (resp.~$C_r^*(\Sigma)$) associated to $\Sigma$, and similarly with $C^*(\B)$ (resp.~$C_r^*(\B)$). For $g\in G$ and $a\in A$, we will identify $\widehat j_g^\A((a,g))$ with the canonical element $a\,U_\Sigma(g)$ in $C^*(\Sigma)$, and $\lambda_g^\A((a,g))$ with the canonical element $a\,\lambda_\Sigma(g)$ in $C_r^*(\Sigma)$.

Consider a family $(\widetilde T_g)_{g\in G}$ of bounded linear maps from $A$ to $B$. For each $g\in G$, define a bounded linear map $T_g:A_g\to B_{\varphi(g)}$ by 
\[ T_g((a, g)) = \big(\widetilde T_g(a), \varphi(g))\]
for all $a\in A$. Then $T=(T_g)_{g\in G}$ is easily seen to be an $\A$-$\varphi$-$\B$ bundle map. Moreover, 
a straightforward computation gives that $T$ is positive definite if and only if the family $(\widetilde T_g)_{g\in G}$ is 
\emph{$\Sigma$-$\varphi$-$\Omega$
positive definite} in the following sense: for every $g_1, \ldots, g_n \in G$ and $a_1, \ldots, a_n \in A$, the matrix
\[  \Big[\beta_{\varphi(g_i)}\big(\widetilde T_{g_i^{-1}g_j}(\alpha_{{g_i}^{-1}}(a_i^*a_j))\big)\Big]\]
is positive in $M_n(B)$.\footnote{When $G=H$, $\varphi={\rm id}_G$, $\Sigma=\Omega$ and $A$ is unital, this amounts to say that $(\tilde T_g)_{g\in G}$ is positive definite w.r.t.~$\Sigma$ in the sense of \cite{BeCo6}.}
Applying Theorem \ref{cp-full-bundle}, we get that $(\widetilde T_g)_{g\in G}$ is $\Sigma$-$\varphi$-$\Omega$ positive definite
 if and only if there exists a completely positive map $\Phi_{\widetilde T}: C^*(\Sigma) \to C^*(\Omega)$ 
such that  \[ 
\Phi_{\widetilde T}(a U_\Sigma(g)) = \widetilde T_g(a) U_\Omega(\varphi(g)) 
\] for all $a\in A$ and $g\in G$, in which case we have $\|\Phi_{\widetilde T}\| = \|\widetilde T_e\|$.
  
If ${\rm ker}(\varphi)$ is amenable, then we may apply Theorem \ref{cp-reduced} and obtain a similar result in the reduced case.
When $G=H$, $\varphi= {\rm id}_G$,  $\Sigma=\Omega$ and $A$ is unital, these two results boil down to the equivalence of (b), (c) and (d) in \cite[Corollary 4.4]{BeCo6}. 
Under the same assumptions, we also illustrated in \cite[Example 4.1]{BeCo6} how positive definite families $(\widetilde T_g)_{g\in G}$ may be constructed from equivariant representations of $\Sigma$. We will now describe how such families arise when $\Omega$ is possibly different from $\Sigma$ and $A$ may be non-unital, still assuming that $G=H$ and $\varphi = {\rm id}_G$.

Let $X$ be a right Hilbert $B$-module such that $A$ acts on $X$ from the left by adjointable operators, and assume $\gamma$ is a homomorphism from $G$ into the group of $\Complessi$-linear invertible isometries of $X$ satisfying
\begin{itemize}
\item[(i)] $\gamma_g(a\cdot x) = \alpha_g(a)\cdot\gamma_g(x)$
\item[(ii)] $\gamma_g(x\cdot b) = \gamma_g(x)\cdot \beta_g(b)$
\item[(iii)] $\langle\gamma_g(x),\gamma_g(y)\rangle_B = \beta_g(\langle x, y\rangle_B)$
\end{itemize}
for all $g\in G, a\in A, b\in B$ and $x, y \in X$. (This means that $\gamma$ is an $\alpha$-$\beta$ compatible action of $G$ on
the right Hilbert $A$-$B$-bimodule $X$ in the sense of \cite{EKQR06}, except that we do not assume non-degeneracy of the left action of $A$ on $X$.) 
For each $g\in G$ define a linear map
$\widetilde T_g : A \to B$ by \[\widetilde T_g(a) = \big\langle x, a\cdot\gamma_g (x)\big\rangle_B \]
for all $a \in A$.
Then one easily verifies that $(\widetilde T_g)_{g\in G}$ is $\Sigma$-${\rm id}_G$-$\Omega$ positive-definite. 
By adapting the proof of \cite[Theorem 4.5]{BeCo6}, it can be shown that every $\Sigma$-${\rm id}_G$-$\Omega$ positive definite family may be obtained in such a way, at least when $A$ and $B$ are unital.

\end{example}

\begin{example} 
 Assume $\B = (B_g)_{g\in G}$ is a Fell sub-bundle of a Fell Bundle $\A=(A_g)_{g\in G}$, as defined in \cite[Definition 21.5]{Exel2}, and $E= (E_g)_{g\in G}$ is a conditional expectation from $\A$ to $\B$ in the sense of Exel \cite[Definition 21.19]{Exel2}. Thus, each $E_g$ is a bounded idempotent linear mapping from $A_g$ onto $B_g$, and we have
\begin{itemize} 
\item[(i)] $E_e$ is a conditional expectation from $A_e$ onto $B_e$, 
\item[(ii)] for each $g, h \in G$ and $a \in A_g, b \in B_h$, we have
\[E_g(a)^* = E_{g^{-1}}(a^*), \quad  E_{gh}(ab) = E_g(a)\,b, \quad E_{hg}(ba)= b E_g(a).\]
\end{itemize}
 Then $E$ is a positive definite $\A$-$\B$ bundle map. 
  Indeed, if $g_1, \ldots, g_n \in G$,  $a_i \in A_{g_i}$ and $b_i\in B_{g_i}$ for $i=1, \ldots, n$, and we set  $c:= \sum_{j=1}^n a_jb_j^* \in A_e$, then 
  \begin{equation*}
 \sum_{i, j=1}^n b_i \, E_{g_i^{-1}g_j} (a_i^* a_j) \, b_j^* = \sum_{i, j=1}^n \, E_{g_ig_i^{-1}g_jg_j^{-1}} ( b_i a_i^* a_jb_j^* ) 
=\sum_{i, j=1}^n \, E_{e} ( b_i a_i^* a_jb_j^* ) = E_e (c^*c) 
\end{equation*}
which is positive in $B_e$ since $E_e$ is a positive map. 

Applying again Theorem \ref{cp-full-bundle} and Theorem \ref{cp-reduced},  we obtain the existence of  completely positive linear maps $\Phi_E: C^*(\A)\to C^*(\B)$ and $M_E:C_r^*(\A) \to C_r^*(\B)$ satisfying that 
\[ \Phi_E(\widehat j^{\,\A}_g(a)) = \widehat j^{\,\B}_g(E_g(a)), \quad M_E(\lambda^\A_g(a)) = \lambda^\B_g(E_g(a))\]
for all $g \in G$ and $a \in A_g$. It is well-known that such a map $M_E$ exists and is a conditional expectation when $C_r^*(\B)$ is identified with its canonical copy inside $C_r^*(\A)$ (cf.~\cite[Theorem 21.29]{Exel2}). In the present situation, we may also  identify $C^*(\B)$ with its canonical copy inside $C^*(\A)$
 (cf.~\cite[Theorem 21.30]{Exel2}), and it is not difficult to check that  $\Phi_E$ is a conditional expectation as well.

\end{example}

\begin{example}\label{C*-corr} 
Let $Y$ be a right Hilbert $C_r^*(\B)$-module and $y\in Y$. Assume that $C_r^*(\A)$ acts on $Y$ from the left by adjointable operators.  
As is well-known and easy to check, we then get a completely positive linear map $M: C_r^*(\A)\to C_r^*(\B)$ by setting
\[ M(X) = \big\langle y , X\cdot y\big\rangle_{C_r^*(\B)} \quad \text{for all } X \in C_r^*(\A).\]
By Proposition \ref{cpred-pdmap}, we get a positive definite $\A$-$\varphi$-$\B$ bundle map  $T=(T_g)_{g\in G}$ by setting
\begin{equation}\label{red-Y-cor} 
T_g(a) := E^\B_{\varphi(g)}\Big( M\big(\lambda^\A_g(a)\big)\Big) = E^\B_{\varphi(g)}\Big(\big\langle y , \lambda^\A_g(a) \cdot y\big\rangle_{C_r^*(\B)}\Big) \quad \text{for all } a\in A_g. 
\end{equation}
Replacing the reduced cross-sectional $C^*$-algebras by their full counterparts and using Proposition \ref{cpfull-pdmap}, we get a similar result by defining $T_g:A_g \to B_{\varphi(g)}$ by 
\begin{equation}\label{full-Y-cor} 
T_g(a) = \big(E^\B_{\varphi(g)}\circ \Lambda^\B\big)\Big(\big\langle y , \widehat j^\A_g(a)\cdot y\big\rangle_{C^*(\B)}\Big) \quad \text{for all } a\in A_g. 
\end{equation}

Conversely, let us start with a positive definite $\A$-$\varphi$-$\B$ bundle map $T=(T_g)_{g\in G}$. 
We can then form the completely positive map $\Phi_T:C^*(\A)\to C^*(\B)$.
Following Paschke \cite[Section 5]{Pas} (see also \cite[p.~48]{La1}),
we may 
use $\Phi_T$ to construct
a right Hilbert $C^*(\B)$-module $Y$ and a left action of $C^*(\A)$ on $Y$ by adjointable operators.
Assume that $A_e$ and $B_e$ are unital. Then $C^*(\A)$ and $C^*(\B)$ are unital, 
so \cite[Theorem 5.2]{Pas} gives that there exists some $y\in Y$ such that
\[\Phi_T(X) = \langle  y, \, X\cdot y\rangle_{C^*(\B)} \]
for all $X\in C^*(\A)$. 
 Then for all $g\in G$ and $a\in A_g$, we have
\[ \langle  y, \, \widehat j^{\,\A}_g(a)\cdot y\rangle_{C^*(\B)} = \Phi_T(\widehat j^{\,\A}_g(a)) = \widehat j_{\varphi(g)}^{\,\B}(T_g(a)),\]
hence 
\[\big(E^\B_{\varphi(g)}\circ\Lambda^\B\big)\big( \langle y, \widehat j^{\,\A}_g(a) \cdot y\rangle_{C^*(\B)}\big)
= E^\B_{\varphi(g)}\big(\Lambda^\B\big(\widehat j_{\varphi(g)}^{\,\B}(T_g(a))\big)\big) 
=   E^\B_{\varphi(g)}\big(\lambda^\B_{\varphi(g)}(T_g(a))\big) = T_g(a) .\]
This shows that when $\A$ and $\B$ are unital, every positive definite $\A$-$\varphi$-$\B$ bundle map $T$ may be written in the form (\ref{full-Y-cor}). 
If we further assume that ${\rm ker}(\varphi)$ is amenable, then arguing analogously with   
$M_T:C_r^*(\A)\to C_r^*(\B)$, we also get that every 
such 
bundle map $T$ is of the form (\ref{red-Y-cor}).
\end{example}

\section{An approximation property for Fell bundles} 
 
 Recall that a discrete group $G$ is amenable if and only if there exists a net $\{\varphi_i\}$ of normalized finitely supported positive definite functions on $G$ such that $\varphi_i \to 1$ pointwise (see e.g.\ \cite[Theorem 2.6.8]{BrOz}). An analogous property for Fell bundles is as follows. 
 
 \begin{definition} A Fell Bundle $\B=(B_g) _{g\in G} $ over a  discrete group $G$ has the \emph{PD-approximation property} if there exists a net $\{T^i\}_{i\in I}$ of $\B$-bundle maps satisfying the following properties:
 \begin{itemize}
 \item[(i)] For each $i \in I$, $T^i=(T^i_g)_{g\in G}$ is positive definite  and its support $\{ g \in G : T^i_g\neq 0\}$ is finite.
  \item[(ii)] The net $\{T^i\}_{i\in I}$ is uniformly bounded in the sense that $\sup_{i}\|T^i_e\| < \infty$.
 \item[(iii)] $\lim_i \|T^i_g(b) - b\| = 0$ for every $g \in G$ and $b\in B_g$. 
 \end{itemize} 
\end{definition}
\begin{example}\label{AP-BC}
Recall from \cite{Exel1, Exel2} that a Fell bundle $\B=(B_g) _{g\in G}$ is said to have  the  \emph{Exel approximation property} (AP)  whenever there exists a net $\{\xi_i\}_{i\in I}$ of finitely supported functions from $G$ to $B_e$ satisfying that 
\begin{itemize}
\item $\sup_{i\in I} \big\| \sum_{h\in G}\xi_i(h)^*\xi_i(h) \big\| < \infty$,
\item $\lim_i \|\sum_{h\in G} \xi_i(gh)^*b\, \xi_i(h) - b\| = 0$ for every $g \in G$ and $b\in B_g$.
\end{itemize}
 It has recently been shown that the AP is equivalent to several other properties, such as AD-amenability, see~\cite{ABF2, BEW2, BFS}. It can be also extended to Fell bundles over groupoids, cf.~\cite{Kra}. 

 We note that a Fell bundle $\B=(B_g)_{g\in G}$ has the PD-approximation property whenever $\B$ has the AP.  Indeed, assume that there exists a net $\{\xi_i\}_{i\in I}$ as above. For each $i\in I$ and $g\in G$, let then $T^i_g: B_g\to B_g$ be the map defined by
 \[ T^i_g (b) = \sum_{h\in G} \xi_i(gh)^* b \, \xi_i(h) \quad \text{for all } b\in B_g\]
As explained in Example \ref{Exel-mult}, each $T^i = (T^i_g)_{g\in G}$ is a positive definite $\B$-bundle map, and one readily
sees that $\{T^i\}_{i \in I}$ is a net guaranteeing that $\B$ has the PD-approximation property. 

This implies in particular that $\B=(B_g)_{g\in G}$ has the PD-approximation property whenever $G$ is amenable. 
\end{example}

\begin{example} 
Let 
$\Sigma=(A, G, \alpha)$ be a unital discrete $C^*$-dynamical system, 
and let $\B$ denote the associated Fell bundle over $G$.  If the system $\Sigma$ is amenable in the sense of \cite{BeCo6}, then, making use of Example \ref{dynsyst}, we get that $\B$ has the PD-approximation property.
\end{example}

 Recall that a Fell bundle $\B$ is said to have the \emph{weak containment property} (WCP), or to be \emph{Exel-amenable}, 
 whenever the canonical map $\Lambda^\B: C^*(\B) \to C_{\rm r}^*(\B)$ is injective. 
\begin{theorem}\label{amen}
Assume  $\B=(B_g) _{g\in G}$ has the PD-approximation property. Then $\B$ has the WCP. Hence $C^*(\B)\simeq C_r^*(\B)$.
\end{theorem}
\begin{proof}
The proof is an adaptation of the proof for the case of group $C^*$-algebras associated to amenable groups given in \cite[Theorem 2.6.8]{BrOz} (see also the proof of \cite[Theorem 4.6] {BeCo6}).

 We start by observing that if $F$ is a finite subset of $G$ and $C_F$ is the subspace of $C_c(\B)$ given by
$C_F= \{ f \in C_c(\B): {\rm supp}(f) \subseteq F\}$, then it is not difficult to check that $\kappa(C_F)$ is a closed subspace  of 
 $C^*(\B)$. Here we write $\kappa$ for the canonical map $\kappa^\B:C_c(\B) \to C^*(\B)$.

Next, let $x \in C^*(\B)$. Consider a positive definite $\B$-bundle map $T$ having finite support. Then $\Phi_T(x) \in \kappa(C_c(\B))$: indeed, letting $F$ denote the support of $T$ and $\{f_n\}$ be a sequence in $C_c(\B)$ such that $\{\kappa(f_n)\}$ converges to $x$, we have (using the notation introduced in Section \ref{prem}) 
\[\Phi_T(\kappa(f_n)) = \kappa(\phi_T(f_n)) =  \kappa\Big(\sum_{g\in F}\big(T_g\big(f_n(g)\big)\odot g\Big)\, \in \kappa(C_F)\] for each $n$, so $\Phi_T(x) = \lim_n \Phi_T(\kappa(f_n))$ belongs to $ \overline{\kappa(C_F)}=\kappa(C_F) \subseteq \kappa(C_c(\B))$.

 Assume now that $\Lambda(x)=0$, where $\Lambda:=\Lambda^\B$ is the canonical $*$-homomorphism from $C^*(\B)$ onto  $C^*_r(\B)$. Let  $\{T^i\}_{i\in I}$ be a net as guaranteed by the PD-approximation property of $\B$. Theorem \ref{cp-full-bundle} (resp.~Theorem \ref{cp-reduced}) gives a net $\{ \Phi_{T^i}\}_{i\in I}$ (resp.~$\{ M_{T^i}\}_{i\in I}$) of completely positive maps on $C^*(\B)$ (resp.\  $C_r^*(\B)$). 
  Note that for each $i\in I$, we have $\Lambda \circ \Phi_{T^i} = M_{T^i}\circ \Lambda$ on $C^*(\B)$, since the maps on both sides are continuous and agree on $\kappa(C_c(\B))$. Hence, for each $i\in I$, we get 
 \[\Lambda \big(\Phi_{T^i}(x)\big) = M_{T^i}\big(\Lambda(x)\big) = 0\,.\]
Since each $T^i$ is finitely supported, we have $\Phi_{T^i}(x) \in \kappa(C_c(\B))$, as established above. But $\Lambda$  is injective on $\kappa(C_c(\B))$, cf.~\cite[Proposition 17.9]{Exel2}, so we obtain that  $\Phi_{T^i}(x) = 0$ for each $i\in I$.

 In order to conclude that $x=0$, we will use the observation that for every $y \in C^*(\B)$ we have $y=\lim_i \Phi_{T^i}(y)$. To see this, let $f\in C_c(\B)$ have support $F\subseteq G$.  Using (\ref{norm-kappa}), we get that
 \[\|\Phi_{T^i}(\kappa(f)) - \kappa(f)\|_{\rm u} = \Big\|\sum_{g\in F} \,\kappa\Big(\big[T_g^i(f(g)) -f(g)\big]\odot g\Big)\Big\|_{\rm u} \, \leq\, \sum_{g\in F} \,\big\|T_g^i(f(g)) -f(g)\|\,.\]
 Hence, using the assumption that $\lim_i \|T^i_g(b) - b\| = 0$ for every $g \in G$ and $b\in B_g$,  we get that $\lim_i \|\Phi_{T^i}(\kappa(f)) - \kappa(f)\|_{\rm u}=0$. This implies that  $\lim_i \|\Phi_{T^i}(y) - y\|_{\rm u}=0$ for all $y \in \kappa\big(C_c(\B)\big)$. If $\{u_j\}_{j\in J}$ is an approximate unit for $B_e$, then $\{\kappa\big(u_j\odot e)\}_{j\in J}$ is easily seen to be an approximate unit for $C^*(\B)$, so we get that  
 \[\|\Phi_{T^i}\| = \lim_j \|\Phi_{T^i}\big(\kappa(u_j\odot e)\big)\|_{\rm u}= \lim_j \|\kappa\big(T_e^i(u_j)\odot e\big)\|_{\rm u} = \lim_j \|T_e^i(u_j)\|  = \|T_e^i\|\]
 for each $i\in I$. 
By property (ii) for $\{T^i\}_{i\in I}$, we have $ \sup_{i\in I} \|T^i_e\| < \infty$, so we get that $\sup_{i\in I}  \|\Phi_{T^i}\| < \infty$. Hence, since $\kappa\big(C_c(\B)\big)$ is dense in $C^*(\B)$, an $\varepsilon/3$-argument gives that
$\lim_i \Phi_{T^i}(y) =y\,$ for every $y \in C^*(\B)$. 

Applying this observation to $x$, we get that 
$ x = \lim_i\Phi_{T^i}(x) = 0$. This shows that $\Lambda $ is injective, as desired.  
\end{proof} 

Our next result involves the tensor product Fell bundle $C \otimes_{\rm min} \B= \big(C \otimes_{\rm min} B_g\big)_{g\in G}$, which is defined for any $C^*$-algebra $C$ in \cite[Definition 25.4]{Exel2}. A similar result when $\B$ has the AP is shown by Exel in \cite[Proposition 25.9]{Exel2}.
\begin{proposition}\label{tensor-BC}
Assume $\B = (B_g)_{g\in G}$ is a Fell bundle having the PD-approximation property, and let $C$ be a $C^*$-algebra. Then  $C \otimes_{\rm min} \B$ has the PD-approximation property.
\end{proposition}
\begin{proof} We recall from \cite[Sections 24 and 25]{Exel2} that  the Fell bundle $\mathcal{D}:=C \otimes_{\rm min} \B$
is the completion of the pre-Fell-bundle obtained by  equipping the algebraic Fell bundle $C\odot \B = (C\odot B_g)_{g\in G}$ with the minimal $C^*$-norm $\|\cdot\|_{\rm min}$ on $C\odot B_e$. The norm on each fiber $C\odot B_g$ is then given by $\|d\| = \|d^*d\|_{\rm min}^{1/2}$ for $d \in C\odot B_g$, and the fiber of $C \otimes_{\rm min} \B$ at $g$, denoted by $C \otimes_{\rm min} B_g$, is the completion of  $C\odot B_g$ w.r.t.~this norm. 

Let  $\{T^i\}_{i \in I}$ be a net as guaranteed by the PD-approximation property of $\B$ and consider $i\in I$. Then by Theorem \ref {cp-reduced}, $T^i$ is a  reduced $\B$-bundle map such that $M_{T^i}$ is completely positive on $C_r^*(\B)$. Letting ${\rm id}_C$ denote the identity map on $C$ and using \cite[Theorem 3.5.3]{BrOz}, we get a completely positive map $\widetilde M_i := {\rm id}_C\otimes M_{T^i}$ on $C \otimes_{\rm min} C_r^*(\B)$ satisfying 
\[\widetilde M_i(c\otimes \lambda^\B_g(b)) = c \otimes M_{T^i}(\lambda^\B_g(b)) =  c \otimes \lambda^\B_g(T^i_g(b))\]
for all $c \in C, g\in G$ and $b \in B_g$, and $\|\widetilde M_i\| = \|{\rm id}_C\| \|M_{T^i}\| = \|T^i_ e\|$.
Now, using \cite[Theorem 25.8]{Exel2}, there is a $*$-isomorphism $\psi: C_r^*(\mathcal{D}) \to C \otimes_{\rm min} C_r^*(\B)$ satisfying  \[\psi ( \lambda^{\mathcal{D}}_g(c\otimes b)) = c \otimes \lambda^\B_g(b)\] for all $c \in C, g\in G$ and $b \in B_g$. 

We therefore get a completely positive map $\Psi_i:= \psi^{-1} \circ \widetilde{M}_i \circ \psi$ on $C_r^*(\mathcal{D})$ satisfying that
\[\Psi_i ( \lambda^\mathcal{D}_g(c\otimes b)) = \lambda^\mathcal{D}_g(c\otimes T^i_g(b))\] for all $c \in C, g\in G$ and $b \in B_g$, and $\|\Psi_i\| = \| T^i_e\|$. 
Hence, using Proposition \ref{cpred-pdmap}, we get that $S^i :=(S^i_g)_{g\in G}$, where for each $g\in G$, $S^i_g: C \otimes_{\rm min} B_g \to C \otimes_{\rm min} B_g$ is defined by 
\[ S^i_g := E^\mathcal{D}_g \circ \Psi_i \circ \lambda^\mathcal{D}_g,\] 
is a positive definite $\mathcal{D}$-bundle map. 

It is now easy to check that the net $\{S^i\}_{i \in I}$ satisfies properties (i) and (ii) needed for the PD-approximation property of  $\mathcal{D}= C \otimes_{\rm min} \B$. To verify that (iii) also holds, using that this net is uniformly bounded, we only have to check that for $c \in C, g\in G$ and $b \in B_g$, we have
\[ \lim_i \|S^i_g(c \otimes b) - c\otimes b\| = 0. \]
Now, as $S^i_g(c \otimes b) = E^\mathcal{D}_g( \Psi_i( \lambda^\mathcal{D}_g(c \otimes b)) = E^\mathcal{D}_g(\lambda^\mathcal{D}_g(c \otimes T^i_g(b))) = c \otimes T^i_g(b)$, we get that 
\begin{align*} 
\|S^i_g(c \otimes b) - c\otimes b\| &= \| (c \otimes (T^i_g(b)-b))^*(c \otimes (T^i_g(b)-b))\|_{\rm min}^{1/2}\\
& =   \| c^*c \,\otimes (T^i_g(b)-b)^*(T^i_g(b)-b)\|_{\rm min}^{1/2} \\ &=  \| c^*c \|^{1/2} \,\|(T^i_g(b)-b)^*(T^i_g(b)-b)\|^{1/2}
= \|c\| \,\|T^i_g(b)-b\| \to_i 0,
\end{align*}
as desired.
\end{proof}

\begin{theorem}\label{amen2}
Assume that $\B=(B_g)_{g\in G}$ is a Fell bundle having the $PD$-approximation property. Then
$C^*(\B) \simeq C_{\rm r}^*(\B)$ is nuclear if and only if $B_e$ is nuclear.
\end{theorem}
\begin{proof}
By Theorem \ref{amen}, we know that $C^*(\B) \simeq C_{\rm r}^*(\B)$. Assume that $B_e$ is nuclear. In \cite[Proposition 25.10]{Exel2}, Exel shows that $C^*(\B) \simeq C_{\rm r}^*(\B)$ is then nuclear whenever $\B$ has the AP. 
It is straightforward to see that his proof goes through verbatim if one just invokes  Proposition \ref{tensor-BC} instead of \cite[Proposition 25.9]{Exel2}.
Conversely, nuclearity of $B_e$ is necessary for the nuclearity of $C_r^*(\B)$ since there exists a conditional expectation from $C_r^*(\B)$ onto $B_e$. 
\end{proof}
We can now include the PD-approximation property 
in \cite[Proposition 7.2]{ABF2} (see also \cite[Th\'eor\`eme 4.5]{AD1}, \cite[Proposition 25.10]{Exel2} and \cite[Corollary 4.23]{BEW2}).  
 \begin{corollary} Let  $\B=(B_g)_{g\in G}$ be a Fell bundle.
  Then the following conditions are equivalent:
 \begin{itemize}
  \item[(i)]$\B$ has the PD-approximation property and $B_e$ is nuclear.
  \item[(ii)]$\B$ has the AP and $B_e$ is nuclear.
 \item[(iii)]$C^*(\B)$ is nuclear.
   \item[(iv)]$C_r^*(\B)$ is nuclear.
  \end{itemize}
 \end{corollary}
 \begin{proof} For the equivalence between ($ii$), ($iii$) and ($iv$), see  \cite[Proposition 7.2]{ABF2}. The implication ($ii$) $\Rightarrow$ ($i$) is shown in Example \ref{AP-BC}. Finally, 
 ($i$) $\Rightarrow$ ($iii$) follows from Theorem \ref{amen2}. 
 \end{proof}
 
 An interesting open question related to this result is whether the PD-approximation property for a Fell bundle $\B=(B_g)_{g\in G}$ is equivalent to the AP, also in the case where $B_e$ is not nuclear.

\section{Some further developments}
In this section we briefly address some directions for further developments.

\subsection{Completely bounded maps on cross-sectional C*-algebras}
Let $\B=(B_g)_{g\in G}$ be a Fell bundle over a discrete group $G$, and let $T$ be a $\B$-bundle map. It seems interesting to look for sufficient conditions 
ensuring that $T$ is  full or reduced with $\Phi_T$ or $M_T$ completely bounded, cf.~Remark \ref{FS} for a result of this type. Here we present another one.

 Given $\phi \in \ell^\infty(G)$,
let $T^\phi$ be the 
$\B$-bundle map given by
  $T^\phi = (T^\phi_g)_{g\in G}$  given by 
\[ T^\phi_g(b) = \phi(g) b\]
for every $g\in G$ and $b\in B_g$ (as in Example \ref{pd-function} when $\phi$ is positive definite). 

Let $B(G)$ denote the Fourier-Stieltjes algebra of $G$ \cite{KL}, consisting of the coefficient functions of unitary representations of $G$. 
If $\phi$ belongs to $B(G)$,  then $T^\phi$ is both full and reduced (as $\phi$ is a linear combination of positive definite functions on $G$).  One may also consider coefficient functions of more general representations of $G$ on Hilbert spaces, as initiated by De Canni\`ere and Haagerup \cite{DCHa} in the case of group bundles. (Note that this will
bring something new 
only when $G$ is non-amenable.) 
The following result extends \cite[Proposition 4.2]{BeCo4} (which is itself an extension of \cite[Theorem 2.2]{DCHa}) to the setting of Fell bundles.

\begin{proposition} Let $v$ be a uniformly bounded representation of $G$ into the invertible bounded linear operators on some Hilbert space $\H$, and set $K:= \sup\{\|v(g)\|, g\in G\} < \infty$. Let $\xi_1, \xi_2 \in \H$ and define 
$\phi \in \ell^\infty(G)$
by
\[\phi(g) = \langle \xi_1, v(g) \xi_2\rangle \quad \text{for all } g \in G.\]
Then $T^\phi$ is reduced. 
In fact, the associated map $M_{T^\phi}:C_r^*(\B)\to C_r^*(\B) $ is completely bounded, with $\|M_{T^\phi}\|_{\rm cb} \leq K^2 \, \|\xi_1\| \|\xi_2\|$.
\end{proposition}
\begin{proof} 
It is not difficult to check that there is an invertible adjointable 
operator $W$ on the Hilbert $B_e$-module $\H \otimes \ell^2(\B)$ satisfying
\[ W(\xi \otimes j^\B_h(c)) = v(h) \xi \otimes j^\B_h(c) \quad \text{ for all } \xi \in \H,  h \in G \text{  and } c \in B_h.\]
A straightforward computation gives that
\[ W (I_\H \otimes \lambda^\B_g(b) )W^{-1} =  v(g) \otimes \lambda^\B_g(b)\]
for every $g\in G$ and $b \in B_g$. 
Now, for $\xi \in \H$, let $\theta_\xi: \ell^2(\B) \to \H \otimes \ell^2(\B) $ be the adjointable operator given by $\theta_\xi f = \xi\otimes f$ for all $f\in \ell^2(\B)$, whose adjoint is determined by $\theta_\xi^* (\xi' \otimes f) = \langle \xi, \xi'\rangle \, f$ for $\xi'\in \H$ and $f \in \ell^2(\B)$.  
We may then define $M_\phi: \L_{B_e}(\ell^2(\B)) \to \L_{B_e}(\ell^2(\B))$ by 
\[ M_\phi(x) = \theta_{\xi_1}^*W(I_\H \otimes x)W^{-1} \theta_{\xi_2}.\]
Then $M_\phi$ is completely bounded with $\|M_{\phi}\|_{\rm cb} \leq K^2 \, \|\xi_1\| \|\xi_2\|$. Moreover, we have
\begin{align*}
 M_\phi (\lambda^\B_g(b)) f & = \theta_{\xi_1}^*W(I_\H \otimes \lambda^\B_g(b))W^{-1} \theta_{\xi_2}f
= \theta_{\xi_1}^*(v(g) \otimes \lambda^\B_g(b)) \theta_{\xi_2}f \\
& = \theta_{\xi_1}^*(v(g)\xi_2 \otimes \lambda^\B_g(b) f)
= \langle \xi_1, v(g) \xi_2\rangle\, \lambda^\B_g(b) f = \phi(g) \, \lambda^\B_g(b) f 
\end{align*}
for all $g \in G, b\in B_g$ and $f \in \ell^2(\B)$. Thus, $M_\phi$, restricted to $C_r^*(\B)$, is a completely bounded extension of $\phi_{T^\phi}$. This shows that $T^\phi$ is reduced, with $M_{T^\phi} = M_\phi$.
\end{proof}

\subsection{Group multipliers} 
Let $G$ and $H$ be discrete groups, $\varphi \in {\rm Hom} (G,H)$, and
 $\phi: G \to \Complessi$ be a bounded function.
Let  us say that $\phi$ is 
a \emph{reduced $(G,\varphi,H)$-multiplier} 
 if there exists a bounded linear map $M_{\varphi, \phi} : C^*_r(G) \to C^*_r(H)$ such that \[M_{\varphi, \phi} (\lambda_G(g)) = \phi(g) \lambda_H(\varphi(g)) \ , \] for every $g \in G$. One can also define \emph{full $(G,\varphi,H)$-multipliers} 
in a similar way,
in terms of the existence a suitable bounded linear map $\Phi_{\varphi, \phi} :C^*(G) \to C^*(H)$.
\noindent Of course, when $G=H$ and $\varphi = {\rm id}_G$,  we recover the usual notions of reduced and full multipliers of $G$ (see e.g.~\cite{BeCo2} and references therein); the map $M_{{\rm id}_G, \phi}$ is then denoted $M_\phi$, while $\Phi_{{\rm id}_G, \phi}$ is denoted by $\Phi_\phi$. In this case, it is well-known (and follows from Theorems \ref{cp-reduced} and \ref{cp-full-bundle}) that the function $\phi$ is positive definite if and only if $\phi$ is a reduced (resp.~full) multiplier and $M_\phi$ (resp.~$\Phi_\phi$) is completely positive.

In terms of the group bundles $\A=(\Complessi \times \{g\})_{g\in G}$ and $\B= (\Complessi \times \{h\})_{h\in H}$, 
letting $T^{\varphi, \phi} = ( T^{\varphi, \phi}_g)_{g\in G}$ be the $\A$-$\varphi$-$\B$ bundle map 
given by
\[ T^{\varphi, \phi}_g(z, g) = (\phi(g) z, \varphi(g)) \quad \text{ for all } z \in \Complessi \text{ and } g \in G,\]
it is clear that $\phi$ is a reduced (resp.~full) $(G,\varphi,H)$-multiplier if and only if the bundle map $T^{\varphi, \phi}$ is reduced (resp.~full).
Moreover, under the natural identifications $C_r^*(\A) \simeq C_r^*(G)$ and $C_r^*(\B) \simeq C_r^*(H)$ (resp.~$C^*(\A) \simeq C^*(G)$ and $C^*(\B) \simeq C^*(H)$), the map $M_{T^{\varphi, \phi}}$ corresponds to $M_{\varphi, \phi}$ (resp.~$\Phi_{T^{\varphi, \phi}}$ corresponds to $\Phi_{\varphi, \phi}$) whenever $\phi$ is a reduced (resp.~full) $(G,\varphi,H)$-multiplier. 

One readily checks that $T^{\varphi, \phi} $ is positive definite if and only if 
 $\phi$ is a positive definite function on $G$. 
Hence, Theorem \ref{cp-full-bundle} gives that  $\phi$ is  positive definite if and only $\phi$ is a full $(G,\varphi,H)$-multiplier such that $\Phi_{\varphi, \phi}:C^*(G)\to C^*(H)$ is completely positive. Moreover, if ${\rm ker}(\varphi)$ is amenable, then  Theorem \ref{cp-reduced}, gives that this is also equivalent to $\phi$ being a reduced $(G,\varphi,H)$-multiplier such that $M_{\varphi,\phi}:C_r^*(G)\to C_r^*(H)$ is completely positive.  

The natural problem of determining the space of all reduced (resp.~full) $(G,\varphi,H)$-multipliers
seems to be challenging. 
We hope to return to this point 
in future work.

\subsection{Fell bundles over groupoids and C*-categories}
Recently there has been a lot of interest in $C^*$-algebras associated  
 to Fell bundles, not only over discrete groups but also over locally compact groups, \'etale groupoids and inverse semigroups (see e.g.~\cite{BE, AF, BEM, KM0, ABF1, KM1, KM2, KLS, KM3, Kra, BM, Tay}).  
It seems natural to investigate to which extent our results in this paper can be carried over to 
these settings. 
Although some work will be necessary in order to handle some technicalities,
we expect that most of our results will continue to hold for Fell bundles over locally compact groups.  In the other cases, it is likely that some efforts will be needed to 
formulate the correct versions. For instance, for   \'etale groupoids, one would first have to decide which class(es) of functors between  \'etale groupoids have the property that a suitable generalization of \cite[Theorems 8.C.12 and 8.C.14]{BedlH} is achievable. We refer to \cite{Tay} for some possible candidates.

In another direction, it has been known for quite some time that Fell bundles over (discrete) groupoids and $C^*$-categories are close relatives, although it is hard to find any specific statement in the literature.
Let $\T$ be a (small) $C^*$-category \cite{GLR}. One can naturally associate to $\T$ a Fell bundle $\B_\T$ over the (discrete) pair groupoid ${\rm Ob}(\T) \times {\rm Ob}(\T)$ by setting 
$B_{\rho,\sigma}:=(\rho,\sigma)$,
the Banach space of arrows from the object $\rho$ to the object $\sigma$ in $\T$. The bundle operations follow at once from the corresponding operations in $\T$. 
If $t \in B_{\rho,\sigma}$ and $s \in B_{\sigma,\tau}$ then $ts \in B_{\rho,\tau}$ (product in the bundle) is given by $s \circ t$ (composition in $\T$), 
and $t^* \in B_{\sigma,\rho}$ (adjoint in the bundle) is nothing but $t^* \in (\sigma,\rho)$. 
Moreover, it is well known that to any $C^*$-category $\T$ one can associate a $C^*$-algebra $C^*(\T)$ (cf.~\cite{GLR}), and to any Fell bundle $\B$ over a groupoid one can associate $C^*(\B)$ and $C_r^*(\B)$ (see e.g.~\cite{Kum}). 
A natural guess would be that $C^*(\T) \simeq C^*(\B_\T)$.

Note also that if $\T$ is a monoidal (i.e., tensor) $C^*$-category, then a second operation of composition in the Fell bundle $\B_\T$ arises from the tensor product in $\T$, 
namely if $t \in B_{\rho,\sigma}$ and $t' \in B_{\rho',\sigma'}$ then we get $t \otimes t' \in B_{\rho \otimes \rho', \sigma \otimes \sigma'}$ (we leave the reader the task of spelling out the detailed properties of this composition). Thinking of a monoidal $C^*$-category as a special case of a 2-$C^*$-category with only one object (see e.g.~\cite{Zi}), this might open the way to a parallel notion of $2$-Fell bundle (cf.~\cite{BCLS}).

Over the years, notions of amenability have appeared in the context of $C^*$-categories mainly inspired by subfactor theory, see e.g.~\cite{LR}  for the concept of amenability of an object in a monoidal $C^*$-category.
In  the more recent paper \cite{PV}, many properties originally defined for groups (e.g., amenability, property (T), the Haagerup property and weak amenability) are 
translated into the setting of rigid monoidal $C^*$-categories.
It is therefore natural to ask if the amenability of a rigid monoidal $C^*$-category $\T$ 
is related to some suitable notion of amenability (or approximation property) of the corresponding Fell bundle $\B_\T$.
It should be possible to put forward several properties from groups to Fell bundles, 
and we  
believe that 
the material developed in the main text 
should provide valuable tools for this purpose.

\bigskip

\noindent{\bf Acknowledgements}

\medskip We thank the referee for her/his careful reading and for making a long list of helpful comments and suggestions.
Erik B\'edos is grateful to the Trond Mohn Foundation (TMS) for partial financial support through the project ``Pure Mathematics in Norway'' during his visit to the Sapienza University of Rome in August 2023.
Roberto Conti (RC) visited the University of Oslo (UiO) in June 2023 and in June 2024. Both these stays were partially supported by UiO. RC was also partially supported by Sapienza Università di Roma and INdAM-GNAMPA project {\it Operator algebras and infinite quantum systems} CUP E53C23001670001.

\section{Statements and Declarations}

\noindent {\bf Competing interests} \ The authors have no Conflict of interest to declare that are relevant to the content of this article. There are no financial and non-financial conflicts of interest.

\medskip

\noindent{\bf Data availability} Data sharing is not applicable to this article as no new data have been created or analyzed in this research.

\bigskip
{\parindent=0pt Addresses of the authors:\\

\smallskip Erik B\'edos, Department of Mathematics, University of
Oslo, \\
P.B. 1053 Blindern, N-0316 Oslo, Norway.\\ E-mail: bedos@math.uio.no. \\
ORCID iD 0000-0002-5559-2571\\

\smallskip \noindent
Roberto Conti, Dipartimento SBAI,
Sapienza Universit\`a di Roma \\
Via A. Scarpa 16,
I-00161 Roma, Italy.
\\ E-mail: roberto.conti@sbai.uniroma1.it\\
ORCID iD 0000-0003-3128-6884
\par}

\end{document}